\documentclass[12pt]{amsart}
\usepackage[colorlinks=true,allcolors=blue,bookmarksopen,bookmarksdepth=3,linktocpage=true]{hyperref}

\usepackage{amsmath, amscd, amsthm, amssymb, graphics, xypic, mathrsfs, setspace, fancyhdr, times, bbold, tikz, tikz-cd, eucal, mathtools,  commath} 
\usepackage[margin=1in,marginparwidth=0.8in, marginparsep=0.1in]{geometry} 

\usepackage{marginnote}
\usepackage{indentfirst}
\usepackage{comment}

\usetikzlibrary{patterns, matrix,shapes,arrows,calc,topaths,intersections,positioning,decorations.pathreplacing,decorations.pathmorphing,fit}

\newcommand{\tensor}{\otimes}
\newcommand{\colim}{\operatorname{colim}}

\newcommand{\cA}{{\mathcal A}}
\newcommand{\cB}{{\mathcal B}}
\newcommand{\cC}{{\mathcal C}}
\renewcommand{\cD}{{\mathcal D}}
\newcommand{\cF}{{\mathcal F}}
\newcommand{\cG}{{\mathcal G}}
\newcommand{\cK}{{\mathcal K}}

\newcommand{\cM}{{\mathcal M}}
\newcommand{\cS}{{\mathcal S}}
\newcommand{\cU}{{\mathcal U}}
\newcommand{\cV}{{\mathcal V}}
\newcommand{\cW}{{\mathcal W}}

\newcommand{\wrap}{\mathfrak{W}}

\newcommand{\cop}{{\operatorname{cop}}}

\renewcommand{\del}{\partial}

\newcommand{\dT}{\dot{T}}
\newcommand{\End}{{\operatorname{End}}}

\newcommand{\F}{{\mathcal F}}

\newcommand{\Fun}{{\operatorname{Fun}}}

\renewcommand{\hom}{\operatorname{Hom}}
\newcommand{\Hom}{{\mathcal Hom}}

\newcommand{\id}{{\operatorname{id}}}
\newcommand{\inhom}{\underline{\operatorname{Hom}} }
\newcommand{\Ind}{{\operatorname{Ind}}}

\newcommand{\Loc}{\operatorname{Loc}}
\newcommand{\Mod}{{\mathcal Mod}}
\newcommand{\ms}{{\operatorname{SS}}}
\newcommand{\msh}{{\operatorname{\mu sh}}}

\newcommand{\real}{{\mathbb R}}
\newcommand{\Shv}{\operatorname{Shv}}

\newcommand{\Sp}{{\mathbb S}}

\newcommand{\supp}{{\operatorname{supp}}}

\newcommand{\Z}{{\mathbb Z}}

\theoremstyle{plain}
\newtheorem{thm}{Theorem}[section]

\newtheorem{lem}[thm]{Lemma}
\newtheorem{cor}[thm]{Corollary}
\newtheorem{prop}[thm]{Proposition}

\newtheorem{construction}[thm]{Construction}

\newtheorem*{thm*}{Theorem}
\newtheorem*{problem*}{Problem}
\newtheorem*{answer*}{Answer}

\theoremstyle{definition}
\newtheorem{defn}[thm]{Definition}
\newtheorem{notation}[thm]{Notation}

\theoremstyle{remark}
\newtheorem{rmk}[thm]{Remark}

\newtheorem{ex}[thm]{Example}

\newtheoremstyle{boldlarge}
  {12pt}
  {12pt}
  {\itshape}
  {}
  {\bfseries}
  {.}
  { }
  {}
\theoremstyle{boldlarge}

\newcommand*{\hrlen}{5}
\newcommand*{\hramp}{3}

\tikzset{
asdstyle/.style={blue,thick},
righthairs/.style={postaction={decorate,draw,decoration={border,amplitude=\hramp,segment length=\hrlen,angle=-90,pre=moveto,pre length=\hrlen/2}}},
lefthairs/.style={postaction={decorate,draw,decoration={border,amplitude=\hramp,segment length=\hrlen,angle=90,pre=moveto,pre length=\hrlen/2}}},
righthairsnogap/.style={postaction={decorate,draw,decoration={border,amplitude=\hramp,segment length=\hrlen,angle=-90}}},
lefthairsnogap/.style={postaction={decorate,draw,decoration={border,amplitude=\hramp,segment length=\hrlen,angle=90}}},
graphstyle/.style={thick},
arrowstyle/.style={thick,decorate,decoration={snake,amplitude=1.7,segment length=10pt,post length=.5mm,pre length=0}},
genmapstyle/.style={thick,-stealth'},     
arrhdstyle/.style={thick},
exceptarcstyle/.style={red, ultra thick},
dualquiverstyle/.style={thick,->},
patstyle/.style={pattern color = gray, pattern = north east lines, opacity=0.3}
}

\title{Functors from Lagrangian correspondences for wrapped sheaf categories}

\author{Sheel Ganatra, Christopher Kuo, Wenyuan Li, Haosen Wu}
\date{}

\begin{document}
\maketitle

\begin{abstract}
We consider the category of sheaves on a manifold with microsupport in a singular conic Lagrangian, and more generally microsheaves on the Lagrangian skeleton of a Weinstein manifold with stop, which are equivalent to ind-completions of partially wrapped Fukaya categories. 
Compact objects in these categories are obtained as wrappings of constructible sheaves with perfect stalks; they correspond to objects built from exact Lagrangians in finitely many steps. 

We give a microlocal criterion on constructible (micro)sheaves with perfect (micro)stalks which guarantees that the corresponding wrapped sheaves induce integral kernels that preserve compact objects of these categories. This generalizes a result of Gao on Lagrangian correspondences between fully wrapped Fukaya categories.
\end{abstract}

\tableofcontents

\section{Introduction}
One powerful invariant of symplectic manifolds is the Fukaya category, proposed by Fukaya \cite{F1}, which is constructed from Lagrangian submanifolds, and counts holomorphic disks  with Lagrangian boundary. Kontsevich in his celebrated homological mirror symmetry conjecture \cite{kon94} %
related this to the category of coherent sheaves on mirror algebraic varieties.

It is a fundamental question 
to study the functoriality of Fukaya categories, both to obtain a functorial invariant on symplectic manifolds and to compare with existing known functoriality results for coherent sheaves on the algebraic varieties. This will be the main subject of this note. 

Let us review the functoriality results for coherent sheaves. %
Algebraic geometers have developed and studied Fourier--Mukai theory on how coherent sheaves on the product space $\cK\in D\mathrm{Coh}(X\times Y)$ give rise to functors 
between the derived category of coherent sheaves on each factor. This was discussed in detail by Mukai, Orlov and Huybrechts \cite{Mukai,Orl02,Huyb}. 
Toen \cite{Toen} first studied the question from the perspective of dg-categories, and Ben-Zvi--Nadler--Francis \cite{BNF} studied the question for derived stacks using $\infty$-categories. 
In the nicest case when the schemes $X,Y$ are both smooth and proper over $S$, their result identified perfect sheaf kernels with exact functors between the categories of perfect complexes: 
$$\Fun^{ex}_{\mathrm{Perf}(S)}(\mathrm{Perf}(X), \mathrm{Perf}(Y)) \cong \mathrm{Perf}(X\times_S Y). $$
When the schemes are not necessarily proper, the situation becomes more subtle. 
Ben-Zvi--Nadler--Preygel \cite{BNP} identified integral kernels that give rise to the exact functors between the derived categories as a subcategory of coherent sheaves on the product $X\times_SY$ %
whose supports are proper over $Y$ and are coherent with finite Tor dimensions relative to $X$: $$\Fun^{ex}_{\mathrm{Perf}(S)}(D^\mathbf{}\mathrm{Coh}(X),D^\mathbf{}\mathrm{Coh}(Y))\cong D_{prop/Y}^\mathbf{}\mathrm{Coh}(X\times_S Y/X).$$ 

Now we turn back to symplectic geometry. 
For closed symplectic manifolds, works of Wehrheim--Woodward, Ma'u--Wehrheim--Woodward \cite{WW,MWW} and later Lekili--Lipyanskiy \cite{LL13} %
pioneered the study of functoriality by showing how to associate to a Lagrangian in the product $M^- \times N$ (where we flip the sign of the symplectic form on $M$) a functor between the Fukaya categories under monotonicity hypotheses. 
The work of 
Fukaya \cite{F} 
generalized the results to any closed Lagrangians, and the works of Bottman and Wehrheim have been further developing the analysis on new types of holomorphic curve bubbling \cite{B15,BW} and the chain-level 2-categorical structure on the category of symplectic manifolds and their Fukaya categories \cite{B4,B5,B25}. %

Fukaya categories of compact symplectic manifolds, when well understood, are often smooth and proper categories \cite{PerutzSheridan,GanatraAuto}, where functors can always be identified as colimit-preserving functors on the ind-completions (see Lemma \ref{lem:proper}). However, when we consider non-compact symplectic manifolds, their (partially wrapped) Fukaya categories are no longer necessarily proper. 

In this paper, %
our goal is to study Fourier--Mukai theory for partially wrapped Fukaya categories of stopped Liouville manifolds (or Liouville sectors) \cite{GPS1}, 
which are often 
non-proper categories. These categories were studied in depth by Ganatra--Pardon--Shende \cite{GPS1,GPS2,GPS3}, %
who in particular show %
an equivalence between partially wrapped Fukaya categories of stopped Weinstein manifolds (or Weinstein sectors) and a purely topological model of these categories in terms of compact objects in the categories of microlocal sheaves on the core of the stopped Weinstein manifolds \cite{NS20}, which goes back to Kashiwara--Schapira \cite{KS13} (see Definition \ref{def:sheaf-microsupport} and \ref{def:microsheaf}):
\begin{equation} \label{for: GPS}
    \mathcal{W}(T^*X, \Lambda)^{op} \cong \Shv_\Lambda^c(X), \quad \mathcal{W}(M, \Lambda)^{op} \cong \msh^c_{\mathfrak{c}_{M,\Lambda}}(\mathfrak{c}_{M,\Lambda}),
\end{equation}
where the superscript $c$ stands for compact objects in a large category (see Definition \ref{def:compact}).
By the K\"unneth formula for such categories \cite{GPS2,Gao,KL25,LiNadlerShende},
we know that colimit-preserving functors between microlocal sheaves or ind-completions of wrapped Fukaya categories are equivalent to microlocal sheaves or ind-completions of wrapped Fukaya categories on the product (see Theorem \ref{thm:kernel} and Corollary \ref{cor:proper-kernel}):
\[\msh_{\mathfrak{c}_{M^- \times N}}({\mathfrak{c}_{M^- \times N}}) \cong \Fun^L(\msh_{\mathfrak{c}_{M}}({\mathfrak{c}_{M}}), \msh_{\mathfrak{c}_{N}}({\mathfrak{c}_{N}})).\]
Hence, the question of understanding functors between the small categories of compact objects again reduces to the question on the integral kernels in the product. 

In the special case of fully wrapped Fukaya categories (with no stops), the question on obtaining functors between the small categories has been previously studied by Y.~Gao \cite{Gao}. 

\begin{thm}[{\cite[Theorem 1.2]{Gao}}]
\label{conj1}
Let $K\subseteq M^- \times N$ be an exact Lagrangian cylindrical at infinity (with Maslov data) whose projection to $N$ is proper.
Then geometric composition with ${K}$ induces a functor between the wrapped Fukaya categories $\cW(M)$ and $ \cW(N)$ (after adding immersed exact Lagrangians with bounding cochains).
\end{thm}

Our work aims to generalize the result to partially wrapped Fukaya categories with non-empty stops, from the sheaf theory perspective. 
Comparing to Fukaya categories, microlocal sheaf categories are well-suited to study functoriality questions via Fourier--Mukai theory, since many functoriality results on categories of sheaves can be packaged into the six functor formalism \cite{KS13,LurieHTT,Volpe}. Moreover, recently, we have developed functoriality results of microlocal sheaves on stopped Weinstein manifolds \cite{LiNadlerShende}.

Consider the stable $\infty$-categories of sheaves with Legendrian microsupports $\Shv_\Lambda(X)$. Recently, the result of Y.~Hu \cite{Hu25} {\em classifies} all sheaves in the product manifold with product stop $\Shv_{-\Lambda \times \Sigma}(X \times Y)$ that induce functors between the {\em compact objects} of the sheaf categories on each factor: (See also Appendix \ref{sec: Hu} for a computational proof through six-functor yoga.)

\begin{thm}[{\cite[Theorem 1.5]{Hu25}}]\label{thm:Hu}
Let $\Lambda \subseteq S^*X$ and $\Sigma \subseteq S^*Y$ be subanalytic Legendrians. Then the integral kernel $\cK \in \Shv_{-\Lambda \times \Sigma}(X \times Y)$ induces a functor preserving compact objects 
\[\Shv_\Lambda^c(X) \to \Shv_\Sigma^c(Y)\]
if and only if for any $x \in X$, the restriction $\cK|_{\{x\} \times Y}$ is compact in $\Shv_\Sigma(Y)$.
\end{thm}

However, for the purpose of understanding Lagrangian correspondences between the wrapped Fukaya categories, we will take a different perspective, using the framework of wrapped sheaves \cite{Kuo}. This is because many interesting functors between the wrapped Fukaya categories are given by Lagrangians in the product that are not asymptotic to the stops, which correspond to constructible sheaves (with perfect stalks) not necessarily microsupported in the stop. Such objects and sheaves that are microsupported in the stop are related by what we call the wrapping functor \cite{Kuo} (see Theorems \ref{thm:wrap-adjoint} and \ref{wrap_compare}):
$$\iota_{\Lambda}^*: \Shv(X) \to \Shv_{\Lambda}(X).$$
Similarly, the same works for microlocal sheaves in general Weinstein manifolds \cite{LiNadlerShende}. We call them wrapped sheaves following \cite{Kuo}. Our question is, therefore, which wrapped sheaves induce functors that restrict to functors between the small categories of compact objects. Our main result is:

\begin{thm}[Main theorem, Theorem \ref{thm:main}] \label{thm: main-intro}
    Let $X, Y$ be real analytic manifolds and $\Lambda \subseteq T^*X$, $\Sigma \subseteq T^*Y$ be conic subanalytic Lagrangian subsets. For any constructible sheaf with perfect stalks $\cK \in \Shv_{\mathbb{R}-c}^b(X \times Y)$, 
    suppose $\ms(\cK) \cap (-\Lambda \times \Sigma) \subset 0_{X \times Y}$, the projection $\supp(\cK) \subseteqq X \times Y \rightarrow X$ is proper and $-\ms(\cK)^v \circ \Sigma$    \footnote{Here $(-)^v$ denotes the coordinate swapping map and $\circ$ the convolution product of Definition~\ref{defn: convolution}.} is negatively squeezable to $\Lambda$ with bounded energy (Definition \ref{def: squeezable}). Then the corresponding wrapped sheaf defines an integral transform preserving compact objects
    \[\iota_{-\Lambda \times \Sigma}^*(\cK) \circ (-): \Shv^c_{\Lambda}(X) \to \Shv^c_{\Sigma}(Y).\]
\end{thm}

Similarly to the main theorem of \cite{GPS3}, the doubling argument reduces the general case to the cotangent bundle case. Let $M, N$ be Weinstein manifolds with mostly Legendrian stops $\Lambda^\infty \subseteq \partial^\infty M, \Sigma^\infty \subseteq \partial^\infty N$. We denote by $\Lambda = \mathfrak{c}_{M,\Lambda^\infty}$ and $\Sigma = \mathfrak{c}_{N,\Sigma^\infty}$ the relative Lagrangian skeleton. Thus, we have the following. 

\begin{thm}[Weinstein case, Theorem \ref{main-Weinstein}] \label{main-Weinstein-intro} 
Let $M, N$ be polarizable Weinstein manifolds with mostly Legendrian stops $\Lambda^\infty \subseteq \partial^\infty M$ and $\Sigma^\infty \subseteq \partial^\infty N$. Let $\cK \in \msh^b_{\real-c}(M^- \times N)$ be a microsheaf with perfect stalks and compact conic Lagrangian support. Suppose 
that $\ms(\cK)^v \circ \Sigma$ is negatively squeezable to $\Lambda$ with bounded energy. Then the functor
\[\iota_{-\Lambda \times \Sigma}^*(\cK) \circ (-):\msh_{\Lambda}(M) \rightarrow \msh_{\Sigma}(N)\]
preserves compact objects.
\end{thm}


Here are a few remarks on the main result. First, this result in particular implies Theorem \ref{conj1} of Y.~Gao \cite{Gao} by Lemma \ref{pre-flowon1factor}. Second, it is not entirely clear to us how to use the classification Theorem \ref{thm:Hu} of Y.~Hu \cite{Hu25} to deduce our main result (however, we are able to deduce a special case of our result from that; see Section \ref{sec:special-case}).

We organize our paper as follows: Section 2 discusses preliminaries on categories, functors and bimodules, Section 3 discusses preliminaries on microlocal sheaf theory, Section 4 recalls sheaf quantizations of Hamiltonian isotopies and wrapped sheaves, Section 5 states and proves the main Theorem \ref{thm:main}, and finally Section 6 implements the main theorem by geometrically constructing certain cofinal wrappings in a special case of Theorems \ref{thm:wrap-one-factor} and \ref{thm:main}.

\subsection*{Acknowledgement}
We would like to thank Yuan Gao, Peter Haine, Yuxuan Hu, David Nadler, and Vivek Shende for helpful correspondence. W.~L. is partially supported by the AMS-Simons Travel Grant.
C.~K. is supported by Max Planck Institute for Mathematics in Bonn. S.G. was partly supported by NSF CAREER DMS-2048055 and NSF FRG DMS-2412041 while working on this project.

\subsection*{AI Disclosure} We used Anthropic's Claude Sonnet solely for grammar polishing, notation consistency, typo checking, and simple calculus computations around Example \ref{ex: Liouville-hypersurface} and Remark \ref{rmk: looser-bounded-energy}. 
The authors are solely responsible for the mathematical content and the final text of the paper.

\section{Higher categories}

\subsection{Compactly generated $\infty$-categories}
We work in the setting of $\infty$-categories 
\cite{LurieHTT}.\footnote{The specific model used in \cite{LurieHTT} is that of quasi-categories, i.e., simplicial sets satisfying the inner horn filling condition. However, in this paper, we are agnostic about the choice of model.} %
We introduce the central notion of \textit{compact objects} in the category.
\begin{defn}[{\cite[Section 5.3.4]{LurieHTT}}]\label{def:compact}
Fix an $\infty$-category $\cC$. An object $C\in \cC$ is \textit{compact} if for any filtered index category $I$ and any functor $I \to \cC$ we have an isomorphism $$\colim_I \hom(C, X_i) \rightarrow \hom\big(C, \colim_I X_i\big).$$
\end{defn}
We state a useful lemma below that relates functors preserving compact objects to properties of their right adjoints: 

\begin{lem}[Preservation of compactness {\cite[Proposition 5.5.7.2]{LurieHTT}}]\label{lem:compact}
Let $\mathcal{C}$ and $\mathcal{D}$ be $\infty$-categories admitting filtered colimits. %
Let $F: \mathcal{C} \to \mathcal{D}$ be a functor that admits a right adjoint $G: \mathcal{D} \to \mathcal{C}$. Then the following conditions are equivalent:
\begin{enumerate}
    \item The left adjoint $F$ preserves compact objects. 
    \item The right adjoint $G$ preserves filtered colimits. 
\end{enumerate}
\end{lem}

\begin{defn}[{\cite[Definition 5.5.7.1 \& 5.5.7.5]{LurieHTT}}]
A category $\cC$ is said to be \textit{compactly generated} if there is a small category $\cC_0$ consisting of compact objects in $\cC$ such that $\cC \cong \Ind(\cC_0)$.
Let $\mathrm{Pr}_\omega^L$ be the subcategory of the category of (large) categories $\widehat{\mathrm{Cat}_\infty}$, where the objects are compactly generated categories and morphisms are functors that preserve small colimits and compact objects.
\end{defn}

Compactly generated categories are special instances of presentable categories:

\begin{defn}[{\cite[Definition 5.5.0.1 \& 5.5.3.1]{LurieHTT}}]
An $\infty$-category $\cC$ is presentable if it admits small colimits, and there is some regular cardinal $\kappa$ such that $\Ind_{\kappa}(\cC_0)=\cC$ for some small category $\cC_0$.
Let $\mathrm{Pr}^L$ be the subcategory of the category of (large) categories $\widehat{\mathrm{Cat}_\infty}$, whose objects are presentable categories and morphisms are functors that preserve small colimits.
\end{defn}




We recall the relation between compactly generated categories and small categories. 


\begin{defn}
Let $\mathrm{Cat}_{\infty}^{Rex}$ be the subcategory of the category $\mathrm{Cat}_{\infty}$ of small categories whose objects are small categories admitting finite colimits and whose morphisms are functors that preserve finite colimits.
\end{defn}
\begin{thm}[{\cite[Proposition 5.5.7.8]{LurieHTT}}]
The functor $\Ind:  \mathrm{Cat}_\infty^{Rex} \to \mathrm{Pr}^L_\omega $ is an equivalence, whose inverse is the nerve of the simplicial functor which associates the compactly generated $\cC$ with the full subcategory $\cC^c \subseteq \cC$ spanned by its compact objects.
\end{thm}

\begin{prop}[{\cite[Proposition 4.8.1.15]{HA}}]%
There exists a symmetric monoidal structure $\tensor$  on $\mathrm{Pr}^L$. Moreover, for enriched compactly generated stable categories $\cC$ and $\cD$, the tensor product $\cC \otimes \cD$ is compactly generated by $c\tensor d$ for compact objects $c\in \cC^c, d\in \cD^c$.
\end{prop}

\subsection{Categorical Fourier--Mukai Theory} 
In this subsection, we quickly recall some categorical notions in the setting of stable $\infty$ or $\cA_\infty$ categories, mainly the definitions of (categorical) smoothness, properness, and relations of bimodules, modules and functors. 

As we mentioned earlier, some early work could be found in \cite{Toen} which describes the $k$-linear continuous functors between two dg-categories $\cD_0,\cC_0$ by $(\cD_0,\cC_0)$-bimodules. 
Here, we follow the notion used in \cite[Appendix A]{GPS3}. %

\begin{defn}[Modules]
A left (respectively right) module over a small category $\cC_0$ is a functor from $\cC^{op}_0$ (respectively $\cC_0$) to the $\infty$-category of spectra $\Sp$. %
A $(\cD_0,\cC_0)$-bimodule over small categories $\cC_0, \cD_0$ is a bilinear functor $\cB$ from $\cD_0^{op} \times \cC_0 \to \cS$. 
\end{defn}

\begin{prop}[{\cite[Proposition 5.3.5.10]{LurieHTT} \cite[Remark 4.8.4.9]{HA}}]\label{prop:bimodule}
Let $\cC_0$ and $\cD_0$ be small stable categories with finite colimits. Then there are equivalences 
\[\cB imod(\cC_0, \cD_0) \simeq \Fun^{ex}(\cC_0, \Ind(\cD_0)) \simeq \Fun^L(\Ind(\cC_0), \Ind(\cD_0)).\]
Moreover, there are also equivalences
\[\Fun^L(\Ind(\cC_0), \Ind(\cD_0)) = \Ind(\cC_0^{op}) \otimes \Ind(\cD_0) = \Ind(\cC_0^{op} \otimes \cD_0).\]
\end{prop}

\begin{defn}[properness]
We say that a $\cC_0$-module $\cM$ is proper if $\cM(C)$ is compact for any $C \in \cC_0$. Denote by $\cC_0^\textit{prop}$ the subcategory of proper modules of $\cC_0$. We say a category $\cC_0$ is \textit{proper} if its diagonal bimodule $\cC_\Delta(-,-)$ is a proper bimodule. 
\end{defn}
\begin{defn}[smoothness]
We say that a $\cC_0$-module $\cM$ is perfect if $\cM$ is a compact module. Denote by $\cC_0^\textit{perf}$ the subcategory of perfect modules of $\cC_0$. We say a category $\cC_0$ is \textit{smooth} if its diagonal bimodule $\cC_\Delta(-,-)$ is a perfect bimodule. 
\end{defn}

\begin{lem}\label{lem:proper}
For small categories with finite colimits $\cC_0$ and $\cD_0$, when $\cC_0$ is smooth but non-proper, there is an embedding
$$\Fun^{ex}(\cC_0, \cD_0) \subseteq \cB imod^\textit{perf}(\cC_0, \cD_0),$$
and when $\cC_0$ is proper but non-smooth, there is an embedding
$$\cB imod^\textit{perf}(\cC_0, \cD_0) \subseteq \Fun^{ex}(\cC_0, \cD_0).$$
\end{lem}
\begin{proof}
The result is essentially the same as \cite[Lemma A.8]{GPS3} \cite[Proposition 4.6.4.4 \& 4.6.4.12]{HA}. When $\cC_0$ is smooth, we have
$$\Fun^{ex}(\cC_0, \cD_0) = \Fun^\textit{prop}(\cC_0, \Ind\cD_0) \subseteq \Fun^\textit{perf}(\cC_0, \Ind\cD_0) = \cB imod^\textit{perf}(\cC_0, \cD_0).$$
Actually, in this case, the identity $id_{\cC_0}$ or equivalently the unit $\eta_{\cC_0}$ is a perfect bimodule. For any proper module $F: \cC_0 \to \cD_0$, since $F$ preserves compactness (and thus perfect), we can deduce that $F = (\epsilon_{\cC_0} \boxtimes F) \circ (id_{\cC_0} \boxtimes \eta_{\cC_0})$ is compact.
Similarly, when $\cC_0$ is proper, we have
\[\Fun^{ex}(\cC_0, \cD_0) = \Fun^\textit{prop}(\cC_0, \Ind\cD_0) \supseteq \Fun^\textit{perf}(\cC_0, \Ind\cD_0) = \cB imod^\textit{perf}(\cC_0, \cD_0).\qedhere\]
\end{proof}

\section{Microlocal sheaves}


\subsection{$\infty$-categories of sheaves}

We describe the $\infty$-category of $\infty$-sheaves (often referred to as the $\infty$-category of sheaves). Let $X$ be a topological space. Denote by $Open(X)$ the poset (viewed as an $\infty$-category) whose objects are open subsets of $X$ and morphisms are inclusions of open subsets. Let $\cC$ be a compactly generated stable $\infty$-category (for example, the derived $\infty$-category $\Mod_R$ of $R$-modules for a given $\mathbb{E}_\infty$-ring spectrum $R$).

\begin{defn} 
Let $X$ be a topological space. A presheaf $\cF$ on $X$ valued in $\cC$ is a functor
\[ 
\cF : Open(X)^{op} \to \cC.
\]
A presheaf $\cF$ is an $\infty$-sheaf if it satisfies descent: for every covering family $\mathcal{U}=\{U_i\}_{i \in I}$ of an open subset $U \subset X$, the canonical map
\[
\cF(U) \longrightarrow \lim_{\Delta^{op}} \cF(\check{C}_\bullet(\mathcal{U}))
\]
is an equivalence, where $\check{C}_\bullet(\mathcal{U})$ is the Čech nerve of the cover $\mathcal U$ and $\Delta$ is the simplex category.
\end{defn}
\begin{rmk}
Concretely, the simplex category $\Delta$ has objects being finite sets of natural numbers $[n]=\{0<\cdots<n\}$, and morphisms being order-preserving maps $[n] \to [m]$. The \v{C}ech nerve is defined by
\[
\check{C}_n(\mathcal{U}) \;=\; \coprod_{(i_0,\dots,i_n)} U_{i_0}\cap\cdots\cap U_{i_n},
\]
with the face/degeneracy maps by inclusions. Applying the functor $\cF$ yields a cosimplicial diagram in $\cC$:
\[
[n]\;\mapsto \;\cF(\check{C}_n(\mathcal{U})) = \coprod_{(i_0,\dots,i_n)} \cF(U_{i_0}\cap\cdots\cap U_{i_n})
\]
and $\lim_{\Delta^{op}}$ is the homotopy limit over this cosimplicial diagram in the $\infty$-category $\cC$.
\end{rmk}

\begin{defn}[$\infty$-category of sheaves]
Let $X$ be a topological space and $\cC$ be a presentable $\infty$-category.  
We denote by $\Shv (X;\cC)$
the $\infty$-category of $\cC$-valued $\infty$-sheaves on $X$ (or the $\infty$-category of $\cC$-valued sheaves on $X$). 
\end{defn}

\begin{thm}[{\cite[Definition 6.1.0.2 \& Corollary 7.1.4.4]{LurieHTT}, \cite[Corollary 2.24]{Volpe}}]
Let $X$ be a paracompact Hausdorff space and $\cC$ be a presentable $\infty$-category. Then $\Shv(X; \cC)$ is a presentable $\infty$-category that admits limits and colimits.
\end{thm}
 
\begin{rmk}
Let $X$ be a paracompact Hausdorff space with finite covering dimensions and $R$ be a discrete ring. Then the $\infty$-category of sheaves $\Shv(X; R)$ is equivalent to the $\infty$-category of hypersheaves, defined as the localization of $\Shv(X; R)$ by maps of sheaves that induce isomorphisms on the stalks by \cite[Theorem 7.2.3.6 \& Corollary 7.2.1.12]{LurieHTT}.

Moreover, the $\infty$-category of sheaves $\Shv(X; R)$ is also equivalent to the derived $\infty$-category of sheaves of abelian $R$-modules on $X$ (whose objects are chain complexes of sheaves and morphisms are maps of complexes of sheaves with quasi-isomorphisms inverted) by \cite[Proposition 7.1]{Scholze}. In particular, the homotopy category $Ho(\Shv(X; R))$ is equivalent to the classical derived category of sheaves of abelian $R$-modules on $X$.
\end{rmk}

We recall that sheaves on locally compact Hausdorff spaces admit a six-functor formalism \cite{Volpe}, which was developed in the classical derived categories in \cite{KS13,Spaltenstein}, and in the dg-derived categories in \cite{SchnurerOlaf}. For this purpose, we assume that the category of compact objects $\cC^c$ is a symmetric monoidal category and that $\cC = \Ind(\cC^c)$ is equipped with the induced presentable symmetric monoidal structure under Ind-completion. In particular, the tensor unit $1_\cC \in \cC^c$ is a compact object.

For a continuous map $f: X \to Y$, we define the direct image functor $f_*: \Shv(X) \to \Shv(Y)$ by
\[\Gamma(V, f_*\cF) = \Gamma(f^{-1}V, \cF).\]
where $\Gamma$ stands for sections. The inverse image functor $f^*: \Shv(Y) \to \Shv(X)$ is defined as its right adjoint. We define the proper direct image functor $f_!$ by
\[\Gamma_c(V, f_!\cF) = \Gamma_c(f^{-1}V, \cF),\]
where $\Gamma_c$ stands for compactly supported sections. The exceptional inverse image functor $f^!: \Shv(Y) \to \Shv(X)$ is defined as its right adjoint
\cite[Section 6.3]{Volpe}.

\subsection{Constructible sheaves}
We define constructible sheaves on a smooth or real analytic manifold $X$. %
In the real analytic setting, we will fix an algebro-geometric category $\cC$, say the category of subanalytic sets in $X$ with respect to the real analytic structure. First, we define the Whitney stratification, see \cite[Part I Chapter 1]{GorMac}:%

\begin{defn} A \textit{Whitney stratification} $\cS = \{\cS_\alpha\}_{\alpha \in A}$ of a smooth manifold $X$ by $\cC$-submanifolds, is a decomposition $X = \sqcup_{\alpha \in A} \cS_\alpha$ into disjoint strata $\{\cS_\alpha\}$ indexed by $A$, where $\cS_\alpha$ are locally closed $\cC$-submanifolds and $\cS_\alpha \cap \overline{\cS_\beta} \neq \varnothing$ if and only if $\cS_\alpha \subset \overline{\cS_\beta}$, and any pair of distinct strata $(\cS_\alpha, \cS_\beta)$ satisfies the \textit{Whitney condition B}.

The Whitney condition B means that: if a sequence $\{x_n \in \cS_\alpha\}$ and a sequence $\{y_n \in \cS_\beta\}$ converge to a point $y \in \cS_\beta$, such that in some local chart the secant lines $\overline{x_n y_n}$ converge to a line $l$ and the tangent planes $T_{x_n} \cS_\alpha$ converge to a plane $\tau$, then $l \subset \tau$.
\end{defn}
\begin{rmk}
Traditionally, we require both Whitney condition A and B, but it is a known result of Mather \cite{Mather} that Whitney condition B implies Whitney condition A.
\end{rmk}
\begin{rmk}
A classical counter-example from \cite{GorMac} for the Whitney condition is given by twisting a rectangular strip sharply on one side. 
\end{rmk}


\begin{defn}
Given a fixed Whitney stratification $\cS$ of a smooth manifold $X$, the $\infty$-category of $\cS$-constructible sheaves
\[
\Shv_\cS (X; \cC)
\]
is the full stable $\infty$-subcategory of $\Shv (X; \cC)$ consisting of objects $\cF$ whose restriction along each stratum of $\cS$ is locally constant and proper objects.
\end{defn}


\begin{rmk}
When $R$ is a discrete ring, all the usual triangulated operations (for example, distinguished triangles, $t$-structures, and truncations) on the classical derived category are lifted %
to the level of the stable $\infty$-categories $\Shv (X;R)$ and $\Shv_\cS (X;R)$.
\end{rmk}



\subsection{Microsupport of sheaves}
Let $X$ be a smooth or real analytic manifold, and let $\cC$ be a compactly generated stable $\infty$-category from now on. Now let us define the microsupport of a sheaf:
\begin{defn}[{\cite[Definition 5.1.1]{KS13}}]
Let $\cF\in \Shv(X; \cC)$. The microsupport of the sheaf $\cF$, denoted by $\ms(\cF) \subseteq T^*X$, is the complement of the set of %
points $p \in T^*X$ such that there exists an open set $U\ni p$, such that for all $x_0\in X$ and all $\varphi: X \rightarrow \real$ satisfying $d\varphi(x_0)\in U$ and $\varphi(x_0)=0$, we have \[\Gamma_{\{x\mid\varphi(x)\geq0\}}(\cF)_x := \mathrm{fib}(\cF \to j_*j^*\cF)_x \cong 0,\]
where $j: \{x \mid \varphi(x) < 0\} \hookrightarrow X$ is the inclusion.

Since the microsupport is conic (invariant under the $\real_+$-action on fibers by scaling), we can define the \textit{microsupport at infinity}: 
$$\ms^\infty(\cF) := (\ms(\cF) \setminus 0_X) / \real_{>0}.$$ 
We will work back and forth between $\ms^\infty(\cF)$ and $\ms(\cF)$. 
\end{defn}

\begin{ex}
We follow closely \cite[Proposition 5.3.1]{KS13} where they provided examples of microsupports of sheaves. 
Fix some open $U = (a,b)$, $a < b$, in $\real$ and inclusion $i: U\hookrightarrow \real$, and let $1_U$ be the constant sheaf on $U$ with value $1_\cC$. We compute 
microsupports of $i_*1_U$ and $i_!1_U$ to illustrate the geometric behavior of microsupports. 

First, we consider $i_*1_U$. We know that $\Gamma(V; i_* 1_U) \simeq \Gamma(U \cap V; 1_U)$, computes the cohomology of $V \cap U$, 
and the stalks can be computed as $(i_*1_U)_x=1_\cC$ if $x\in {U}$ and $(i_*1_U)_x=0$ if $x\notin {\overline{U}}$. 
The microstalk can then be computed as $m_{(x,\xi)}(i_*1_U)=1_\cC$ if $x=a,\xi\geq0$ or $a<x<b,\xi=0$ or $x=b,\xi\leq0$. Therefore the microsupport is 
\[\ms^\infty(i_*1_U)=\{(a,\xi)\in S^*\real \mid \xi>0\}\cup\{(b,\xi)\in S^*\real \mid \xi<0\}.\]

Now, let us investigate $i_!1_U$. We know that $\Gamma(V; i_! 1_U) \simeq \Gamma_c(V \cap U; 1_\cC)$ computes the relative cohomology of the pair $(V, V \setminus U)$, and the stalks can be computed as $(i_!1_U)_x=1_\cC$ if $x\in U$ and $(i_!1_U)_x=0$ if $x\not\in U$. 
The microstalk can be computed as $m_{(x,\xi)}(i_!1_U)=1_\cC[-1]$ if $x=a,\ \xi<0$ or $x=b,\ \xi>0$. For $a< x< b$ and $\xi=0$, it is simply the stalk, so $m_{(x,\xi)}(i_!1_U)=1_\cC$.
Therefore  
the microsupport is \[\ms^\infty(i_!1_U)=\{(a,\xi)\in S^*\real\mid \xi<0\}\cup\{(b,\xi)\in S^*\real\mid \xi>0\}.\] 
\end{ex}

\begin{ex}
We give more delicate examples in terms of graphs, following \cite{Zhou}: Let $j: U=B(0,1) \hookrightarrow \real^2$ be the inclusion of an open unit ball in $\real^2$. Then $j_* 1_U$ has microsupport at infinity as 
\[ \ms^\infty(j_* 1_U) = \{ (x, \xi) \in S^*\real^2 \mid |x|=1, \xi = -d(|x|^2) \} = \quad \tikz[baseline=-3pt] {\draw [lefthairs] (0,0) circle (0.5);}. \]
$j_! 1_U$ has microsupport at infinity as
\[ \ms^\infty(j_! 1_U) = \{ (x, \xi) \in S^*\real^2 \mid |x|=1, \xi = d(|x|^2) \} = \quad \tikz[baseline=-3pt] {\draw [righthairs] (0,0) circle (.5);}.\]
\end{ex}
In summary, microsupports document where and in which codirection the (derived) sections of the sheaves change. We refer to such behavior as the failure of propagation of sections of the sheaves. By \cite[Theorem 6.5.4]{KS13}, the microsupport of a sheaf is always coisotropic and, at smooth Lagrangian points, this non-propagation can be measured more accurately using microstalk:

\begin{defn}
    Let $(x,\xi) \in \ms(\cF)$ and assume that $\ms(\cF)$ is a smooth Lagrangian manifold near $(x,\xi)$. The \textit{microstalk} of a sheaf $\cF$ at $(x, \xi)$, denoted as $m_{(x,\xi)}\cF$, with respect to a smooth function $\varphi$ defined near $x$ with $\varphi(x) = 0$ and $d\varphi(x) = \xi$, is defined as \[m_{(x,\xi)}\cF \coloneqq \Gamma_{\{x\mid \varphi(x) \geq 0\}}(\cF)_x.\]
\end{defn}
\begin{rmk}
    We warn the readers that the microstalk of $\cF$  at $(x, \xi) \in T^*X$ is not independent of the choice of the smooth function $\varphi$. However, generically, it is well-defined up to a cohomological degree shift, and one can calculate the shift with respect to different choices of the smooth functions by \cite[Section 7.4]{KS13}.
\end{rmk}



\subsection{Sheaves with prescribed microsupports; Microsheaves}
We introduce the categories of sheaves with prescribed microsupport: 
\begin{defn}\label{def:sheaf-microsupport}
    For a closed conic subset $\Lambda \subseteq T^*X$,  the category of sheaves with prescribed microsupport, $\Shv_\Lambda (X; \cC) $ is the full subcategory of $\Shv(X; \cC)$ consisting of objects with microsupport contained in $\Lambda$. Similarly, for a closed subset $\Lambda \subseteq S^*X$, $\Shv_\Lambda (X; \cC)$ is the full subcategory of $\Shv(X; \cC)$, consisting of objects %
    with microsupport at infinity contained in $\Lambda$.
\end{defn}
\begin{thm}[{\cite[Lemma 4.12]{GPS3} \cite[Proposition 3.4]{GuiVit}}]\label{thm: adjoint}
For a closed conic subset $\Lambda \subseteq T^*X$ (or a closed subset $\Lambda \subseteq S^*X$), the inclusion functor
\[\iota_{\Lambda*}: \Shv_\Lambda(X; \cC) \hookrightarrow \Shv(X; \cC)\]
is limit- and colimit-preserving, and hence admits a left adjoint $\iota_\Lambda^*$ and a right adjoint $\iota_\Lambda^!$. In particular, $\Shv_\Lambda(X; \cC)$ is a stable presentable category.
\end{thm}

We can make a comparison between the category of sheaves with prescribed microsupport and the category of constructible sheaves. Since the result ultimately relies on \cite[Proposition 8.3.10]{KS13}, which states that if $\ms(F)$ is subanalytic, then $F$ is constructible if and only if $\ms(F)$ is a (possibly singular) Lagrangian, we will restrict our attention to the subanalytic subsets. That is, unless specified, any stratification $\cS$, any Lagrangian $\Lambda \subseteq T^* X$, etc. will be subanalytic.

\begin{thm}[{\cite[Theorem 8.4.1]{KS13} \cite[Theorem 4.8]{GPS3}}] \label{thm: ss-constructible}
For a Whitney stratification $\cS$ of a $C^1$ manifold $X$, we have an equivalence ${\Shv}_{\cS}(X; \cC) = {\Shv}_{N^*\cS}(X; \cC)$.
\end{thm}

We introduce the subcategory of compact and proper objects in $\Shv_\Lambda^c(X; \cC)$:
\begin{defn}
Let $\Shv_\Lambda^c(X; \cC)$ be the full subcategory of compact objects in $\Shv_\Lambda(X)$ and let $\Shv_\Lambda^\textit{pr}(X; \cC)$ be the full subcategory of proper modules of $\Shv_\Lambda^c(X; \cC)$. 
\end{defn} 

\begin{thm}[{\cite[Corollary 4.21 \& 4.23]{GPS3}}]\label{thm:smoothness}
Let $X$ be compact and $\Lambda$ be subanalytic Lagrangian. Then $\Shv_\Lambda(X; \cC)$ is a smooth category compactly generated by corepresentatives of microstalks, and $\Shv_\Lambda^\textit{pr}(X; \cC)$ is equivalent to sheaves in $\Shv_\Lambda(X; \cC)$ with perfect stalks.
\end{thm}  

A particularly nice case is when $\Lambda = N^* \cS$ where $\cS$ is a triangulation.
\begin{lem}[{\cite[Lemma 4.18]{GPS3}}] \label{lem: smoothness-triangulation}
 For $\cS$ a triangulation, the category $\Shv_\cS^c(X)$ consists of $\cS$-constructible sheaves with compact support and perfect stalks.
\end{lem}

We can use the categories $\Shv_\Lambda(X)$ characterized by the microsupports to define a sheaf of categories on the cotangent bundle $T^*X$, called the \textit{microsheaf category}:

\begin{defn}[Kashiwara--Schapira stack/the microsheaf category]\label{def:microsheaf} Define the presheaf of categories on $T^*X$ with respect to the conic topology to be such that for $\Omega \subset T^*X$ a conic open subset, 
\[\msh^{pre}(\Omega; \cC):=\Shv(X; \cC)/\Shv_{T^*X\backslash \Omega}(X; \cC).\]
Then we define $\msh$ to be the sheafification of $\msh^{pre}$ valued in the $\infty$-category of all (not necessarily small) stable $\infty$-categories.
    
For a closed conic subset $\Lambda \subseteq T^*X$ (or similarly a closed subset $\Lambda \subseteq S^*X$), 
we define $\msh_{\Lambda}(-)$ to be the subsheaf of $\msh(-)$ whose sections are full subcategories of objects with (micro)supports in $\Lambda \subseteq T^*X$ (or respectively $\Lambda \subseteq S^*X$). 
\end{defn}
\begin{rmk}
We require the closedness of conic subset $\Lambda$ to ensure that the category is closed under limits and colimits, 
which is necessary to take Bousfield localization  \cite{hirschhorn}. 

However, the sheafification process takes place in the category of stable categories and may not preserve limits and colimits. The reason to sheafify in the $\infty$-category of categories is that it admits a limit-preserving functor to the $\infty$-category of spaces \cite{Rozenblyum-filtered}, so that one can check isomorphisms of sheaves on the stalks \cite[Theorem 7.2.3.6 \& Corollary 7.2.1.12]{LurieHTT}. This allows one to identify the morphisms of $\msh$ with the $\mu hom$ functor by \cite[Theorem 6.1.2]{KS13}.
\end{rmk}

\begin{rmk}
When $\Lambda \subseteq T^*X$ is a conic subset containing the zero section
, plugging $\Lambda$ into the sheaf of categories, 
we can show that $\Shv_\Lambda(X; \cC) \cong \msh_\Lambda(T^*X; \cC) \cong \msh_\Lambda(\Lambda; \cC)$. In fact, 
\[\msh_\Lambda(T^*X; \cC) \cong \pi_*\msh_\Lambda(X; \cC) \cong (\pi_*\msh_\Lambda^{pre})^{sh}(X; \cC) \cong \Shv_\Lambda(X; \cC),\]
since the pre-sheaf of categories of sheaves on $X$ forms a sheaf of categories \cite[Theorems 6.1.0.6, 6.1.3.9 \& Proposition 6.2.2.7]{LurieHTT} \cite[Corollary B.5]{Jin-Treumann}. 
\end{rmk}

\begin{prop}[{\cite[Proposition 3.16]{CasalsLi} \cite[Proposition 5.5.3]{LiNadlerShende}}]\label{prop:musheaf-lim-colim}
Let $\Lambda \subseteq S^*X$ be a precompact subanalytic Legendrian. Then $\msh_\Lambda$ is a sheaf valued in $\mathrm{Pr}^L_\omega$, where sections are compactly generated and restriction functors preserve limits and colimits.
\end{prop}

Let $M$ be a Weinstein manifold with subanalytic Lagrangian skeleton $\mathfrak{c}_M$ as in \cite[Section 7.7]{GPS3}, equipped with a polarization $\sigma$ which is a null homotopy of the Lagrangian Gauss map. For $N \gg 0$, we can embed the contactization $M \times \real $ into $S^*\real^N$ as a contact submanifold by the $h$-principle, and the polarization determines a Lagrangian subbundle in the normal bundle \cite[Section 9.1]{NS20}. Denote the restriction of the Lagrangian section to the skeleton in $M \times 0$ by $\mathfrak{c}_{M,\sigma}$. Then by the exponential map, we have a subanalytic Lagrangian $\mathfrak{c}_{M,\sigma} \subseteq S^*\real^N$.

\begin{defn}[{\cite[Definition 9.4]{NS20}}] \label{defn: microsheaf-Weinstein}
Let $M$ be a Weinstein manifold with subanalytic Lagrangian skeleton $\mathfrak{c}_M$. Given a polarization $\sigma$, the microsheaf category of $M$ on the Lagrangian skeleton is
$$\msh_{\mathfrak{c}_M} = i^*\msh_{\mathfrak{c}_{M,\sigma}},$$
where $i: \mathfrak{c}_{M} \hookrightarrow \mathfrak{c}_{M,\sigma}$ is the inclusion of the zero section. 
\end{defn}

We briefly recall the key ingredient used to study microsheaves: the doubling functor. For a Legendrian $\Lambda \subseteq \dot T^*X$, and a Reeb flow $T$ on $\dot T^*X$, recall that the Reeb movie of $\Lambda$ with length $\epsilon$ is
\[\Lambda_{\cup,\epsilon} = T_{-\epsilon}(\Lambda) \cup T_\epsilon(\Lambda) \cup \pi\Big(\bigcup\nolimits_{s\in[-\epsilon,\epsilon]}T_s(\Lambda)\Big) \subset T^*X.\]
Let $\Lambda_{\cup,\epsilon}^\pm$ be the Reeb push-off of $\Lambda_{\cup,\epsilon}$ by $\pm\epsilon$.

\begin{prop}[{\cite[Theorem 5.2.1]{LiNadlerShende}}]\label{prop: doubling}
Let $M$ be polarizable a Weinstein manifold with mostly Legendrian stop $\Lambda^\infty$, and the relative Lagrangian skeleton be $\Lambda$. There for large enough $N\in \mathbb{N}$, the category $\msh_{\Lambda}(M)$ embeds as a presentable subcategory of $\Shv(\real^N)$. More precisely, there is an equivalence 
\[w_{\Lambda}^\pm: \msh_{\Lambda}(M) \xrightarrow{\sim} \Shv_{(\Lambda_{\sigma})_{\cup,\epsilon}^\pm}(\real^N) \subseteq \Shv(\real^N).
\] 
\end{prop}



\subsection{Bimodules of microlocal sheaves}
Colimit preserving functors, or equivalently bimodules by Proposition \ref{prop:bimodule}, on the categories of microlocal sheaves can be classified by integral kernels in the product. We recall the results:

\begin{thm}[{\cite[Theorem 1.1 \& 1.4]{KL25} \cite[Theorem B]{Zhang}}]\label{thm:kernel}
Let $\Lambda \subseteq T^*X$, $\Sigma \subseteq T^*Y$ be closed conic subsets. Then there is an equivalence
$$\Shv_{-\Lambda \times \Sigma}(X \times Y) \cong \Fun^L(\Shv_\Lambda(X),\Shv_\Sigma(Y)) \cong \cB imod(\Shv_\Lambda^c(X),\Shv_\Sigma^c(Y)),$$
that sends $\cK$ to the integral transform $\cK \circ (-) \coloneqq \pi_{Y!}(\cK \otimes \pi_X^*(-))$. Similarly, let $\Lambda \subseteq S^*X$, $\Sigma \subseteq S^*Y$ be closed subanalytic Legendrian subsets. Then there is an equivalence
$$\msh_{-\Lambda \times \Sigma}(-\Lambda \times \Sigma) \cong \Fun^L(\msh_\Lambda(\Lambda),\msh_\Sigma(\Sigma)) \cong \cB imod(\msh_\Lambda(\Lambda),\msh_\Sigma(\Sigma)),$$
that sends $\cK$ to the integral transform.
\end{thm}

\begin{notation} \label{not: circ-hom}
For a sheaf kernel $\cK$, we set
\[ \Hom^\circ(\cK, \cG) = \pi_{X*}\Hom(\cK, \pi_Y^!\cG). \]
Note that $\Hom^\circ(\cK,-)$  is the right adjoint of $\cK \circ (-)$.
\end{notation}

As a corollary of Theorem \ref{thm:smoothness} and Lemma \ref{lem:proper}, we know that all integral kernels that preserve compact objects are necessarily compact objects in the product:

\begin{cor}\label{cor:proper-kernel}
Let $X$ be compact and $\Lambda \subseteq T^*X$, $\Sigma \subseteq T^*Y$ be closed conic subsets. Then there is an embedding
$$\Fun^{ex}(\Shv_\Lambda^c(X),\Shv_\Sigma^c(Y)) \subseteq \Shv_{-\Lambda \times \Sigma}^c(X \times Y).$$
Similarly, let $\Lambda \subseteq S^*X$, $\Sigma \subseteq S^*Y$ be closed subanalytic Legendrian subsets. Then there is an embedding
$$\Fun^{ex}(\msh_\Lambda^c(\Lambda),\msh_\Sigma^c(\Sigma)) \subseteq \msh_{-\Lambda \times \Sigma}^c(-\Lambda \times \Sigma).$$
\end{cor}

\begin{defn}
We say a kernel $\cK \in \Shv_{-\Lambda \times \Sigma}(X \times Y)$ is compact-object-preserving if the associated functor is compact-object preserving. Furthermore, we use the notation $\Shv_{-\Lambda \times \Sigma}(X \times Y)^\cop \subseteq \Shv_{-\Lambda \times \Sigma}(X \times Y)$ for the stable subcategory consisting of compact-object-preserving kernels. 
\end{defn}

We note that without compactness assumption on $X$, $\Shv_{-\Lambda \times \Sigma}(X \times Y)^\cop$ neither contains nor is contained in the category of compact objects $\Shv_{-\Lambda \times \Sigma}(X \times Y)^c$.

\begin{ex} \label{ex: cop-not-c}
Consider the case when $X = Y = \real$ the real line and $\Lambda = \Sigma = N^* \cS$ where $\cS$ is the stratification given by $\cS = \{\Z, \real \setminus \Z\}$. In this case, colimit-preserving endofunctors $\End^L(\Shv_{N^* \cS}(\real)) = \Shv_{\cS \times \cS}(\real \times \real)$ are given by $\cS \times \cS$-constructible sheaves on $\real^2$. Per Lemma \ref{lem: smoothness-triangulation}, compact objects in particular have compact support. However, as $\Lambda$ is swappable in the sense of \cite[Definition 5.19]{KL22}, the wrap-once functor $S_\Lambda^+(F) \coloneqq \wrap_\Lambda^+(F^w)$, where $(-)^w$ is any small wrapping displacing $N^*(\Z)$ positively from itself, is an equivalence so in particular preserves compact objects. In our case, $S_\Lambda^+$ corresponds to $\iota_{N^*(\cS) \times N^*(\cS)}^*(1_{\{|y-x| < \epsilon\}})[1]$, whose support contains the diagonal $\Delta$ as a closed subset. Thus, it is not a compact object. 
\end{ex}

\begin{ex}
Consider the case when $X = S^1$, $Y = *$, and $\Lambda = 0_{S^1}$ be the zero section of $S^1$. Denote by $\pi: \real \rightarrow S^1 = \real/ \Z$ the universal cover. One can check directly that $\pi_! 1_{\real}$ is a compact generator. Its corresponding functor is given by
\begin{align*}
\Loc(S^1) &\rightarrow \cV \\
L &\mapsto \Gamma(S^1; (\pi_! 1_{\real}) \otimes L).
\end{align*}
This last expression is not compact-object-preserving. For example, take $L = \pi_! 1_{\real}$ itself and one can compute, by the projection formula, that 
\[ \Gamma(S^1; (\pi_! 1_{\real}) \otimes \pi_! 1_{\real}) = \Gamma(S^1; \pi_! (\pi^* \pi_! 1_{\real})) = \Gamma(\real; \pi^* \pi_! 1_{\real}) \cong i_0^*  \pi^* \pi_! 1_{\real} = 1_{\real}^{\oplus \Z}   \]
is not a compact object. Here, we use the fact that, on $\Loc(\real)$, taking global section $\Gamma(\real;-) \cong i_0^*$ has the same effect as taking stalk at $0$.

In fact, the same computation implies that a kernel $\cK \in \Loc(S^1)$ is compact-object-preserving if and only if it has perfect stalk. This is because compact objects are closed under finite colimit and retracts, so $\cK$ preserves compact objects as long as it sends $\pi_! 1_{\real}$ to a compact object and 
\[ \Gamma(S^1; \cK \otimes \pi_! 1_{\real}) =  \Gamma(\real; \pi^* \cK) \cong i_0^*  \pi^* \cK \cong \cK_{\bar{0}}.   \]
\end{ex}

\section{Isotopy techniques in microlocal sheaf theory}

\subsection{Guillermou--Kashiwara--Schapira sheaf quantization}
In this section we first review the construction of Guillermou--Kashiwara--Schapira sheaf quantization (or in short, GKS sheaf quantization), which gives rise to the foundation of sheaf-theoretic wrapping, %
and then review the definition of sheaf-theoretic wrapping along with necessary notions, defined in \cite{GKS} and further systematically studied in the second author's thesis \cite{Kuo}. %


\begin{defn}[Contact Hamiltonian flow]\label{positivewrapping}
Let $(S^*X,\alpha)$ be a contact manifold, and let $R_\alpha$ be the Reeb vector field.
A time-dependent contact Hamiltonian (or simply contact isotopy) is a smooth function
$H : S^*X \times I \to \mathbb{R}$, and we write $H_t(y) := H(y,t)$.
For each $t$, the associated contact Hamiltonian vector field $V_{H_t}$ is determined by
\begin{equation} \label{for: contact-flow}
    \alpha(V_{H_t}) = H_t, \qquad
\iota_{V_{H_t}} d\alpha = (dH_t)(R_\alpha)\,\alpha - dH_t,
\end{equation}

and its flow $\phi_H^t$ is the contact isotopy solving
$\partial_t \phi_H^t = V_{H_t} \circ \phi_H^t$ with $\phi_H^0 = \mathrm{id}_{S^*X}$. It is called \emph{positive} if $H_t \ge 0$. Such positive isotopies are called wrappings. 

Let $H: S^*X \times I \to \real$ be a contact Hamiltonian. We abuse notation and write the associated homogeneous Hamiltonian as $H: \dT^*X \times I \to \real$, which is given by the formula 
\[H(y, r, t) = rH(y, t), \quad (y, r) \in S^*X \times \real_{\ge 0} / S^*X \times \{0\}, \; t \in I.\]
Its flow $\phi_H^t$ is the symplectic Hamiltonian isotopy on $\dT^*X$ defined by $H$.\footnote{We will frequently pass between the contact and homogeneous languages using this construction.}
\end{defn}

\begin{rmk}
We note that in the homogeneous setting, there is a unique homogeneous Hamiltonian function. Let $\Phi: \dT^*X \times I \rightarrow \dT^*X$. Its homogeneous Hamiltonian function is given by $H(x,\xi,t) \coloneqq \alpha_{\varphi_t(x,\xi)}(\partial_t \varphi_t(x,\xi))$.
\end{rmk}

\begin{thm}[Guillermou--Kashiwara--Schapira sheaf quantization, {\cite[Theorem 3.7]{GKS}}]\label{thm:gks}
Let $\Phi$ be a homogeneous %
Hamiltonian isotopy of $\dT^*X$ as above. Then there exists a sheaf kernel
$\cK(\Phi)\in \Shv(X\times X\times I)$
with the following properties:
\begin{itemize}
\item[1)] $\ms^\infty(\cK(\Phi))$ is contained in the conical Lagrangian, the ``movie" of $\Phi$,
$$\Lambda_\Phi = \{(x, -\xi, \varphi_t(x, \xi), t, -H(\varphi_t(x, \xi), t) \mid (x, \xi) \in \dT^*X, t\in I\};$$
\item[2)] $\cK(\Phi)$ has microlocal rank $1$ along $\Lambda_\Phi$;
\item[3)] $\cK(\Phi)|_{X\times X \times\{0\}}\simeq 1_{\Delta_M}$.
\end{itemize}
Moreover, $\cK(\Phi)$ is unique up to isomorphism under these conditions. We will call this sheaf kernel $\cK(\Phi)$ the GKS kernel.
\end{thm}

\begin{rmk}[{\cite[Remark 3.9]{GKS}}] \label{rmk: family-GKS}
    In fact, $I$ can be replaced by any contractible manifold $U$ with a distinguished point $u_0$ such that $\Phi(-,u_0) = \id_{\dT^*X}$. In this case, one replaces the movie $\Lambda_\Phi$ by the Lagrangian lift of the Legendrian movie
    \[L_\Phi = \{(x, -\xi, \varphi_u(x, \xi), u) \mid (x, \xi) \in \dT^*X, u\in U\}.\]
    For example, suppose $U = I \times J$ is the product of two intervals with coordinates $t$ and $s$, and denote by $\tau$ and $\sigma$ the corresponding cotangent coordinates. Then $\Lambda_\Phi \subseteq \dT^*X \times \dT^*X \times T^*I \times T^*J$ satisfies
    \begin{equation} \label{for: family-Hamiltonian-function}
    \tau = -\alpha_{\varphi_{t,s}(x,\xi)}(\partial_t\varphi_{t,s}(x,\xi)) \ \text{and} \ \sigma =-\alpha_{\varphi_{t,s}(x,\xi)}(\partial_s\varphi_{t,s}(x,\xi)).    
    \end{equation}
\end{rmk}

These GKS kernels act on sheaves through convolution $\circ$.


\begin{prop}[Effect of the GKS kernel {\cite[Proposition 3.12]{GKS}}]\label{prop:gks-effect}
Let $\Phi$ be a homogeneous Hamiltonian isotopy and $\cK(\Phi)\in \Shv(X\times X\times I)$ be the GKS kernel. Write $\cK(\Phi)|_t := \cK(\Phi)|_{X\times X\times\{t\}}$.
Then:
\begin{itemize}
\item[1)] $\cK(\Phi)|_0 \circ \cF \simeq \cF$;
\item[2)] 
${\ms}^\infty(\cK(\Phi)|_t \circ \cF)\subset {\varphi_t({\ms}^\infty(\cF))}$ for all $t$;
\item[3)] convolution with $\cK(\Phi)|_t$ restricts to an equivalence
$\Shv_\Lambda(X)\xrightarrow{\sim} \Shv_{\varphi_t(\Lambda)}(X)$
with inverse given by convolution with the kernel corresponding to $\varphi_t^{-1}$.
\end{itemize}
\end{prop}

\begin{ex}
Take $X=\real$ and let $\varphi_t$ be the contact isotopy on $S^*\real$ that sends $(x, \xi)$ to $(x-t, \xi)$ if $\xi > 0$ and $(x+t, \xi)$ if $\xi < 0$. Let $\Phi$ denote the corresponding family $\Phi(-,t)=\varphi_t$. Let $\cF = 1_0$ be the skyscraper sheaf whose microsupport at infinity is the fiber of $0 \in \real$. Then 
$$\cK(\Phi)|_t = 1_{\{(x, y) \mid |x-y|\leq t\}}, \quad \cK(\Phi)|_t \circ \cF = \cK_t\circ \cF = 1_{[-t,t]},$$
and $\cK(\Phi)|_t \circ \cF$ has microsupport at infinity along the positive ray in the fiber of $-t \in \real$ and the negative ray in the fiber of $t \in \real$. In this way the GKS kernel transports the microsupport along the contact flow.
\end{ex}

\begin{prop}[{\cite[Proposition 4.8]{GKS} \cite[Proposition 3.2 \& 3.9]{Kuo}}] \label{prop: continuation-map}
Let $\Phi$ be a homogeneous positive Hamiltonian isotopy. Then 
$\ms(\cK(\Phi)\circ \cF)\subset T^*X\times T^*_{\le 0}I$.
There is a canonical continuation morphism
$$\cK(\Phi)|_a \circ \cF \to \cK(\Phi)|_b \circ \cF, \quad a\le b.$$ 
\end{prop}
For simplicity, we set $\cF_a \coloneqq \cK(\Phi)|_a \circ \cF$ when the isotopy is clear from the context. These continuation maps are used in wrapping: as time increases, the sheaf is transported by the isotopy, and the maps $\cF_a\to \cF_b$ record this motion in a way compatible with microsupport and composition in $\Shv(X)$. As explained in \cite[Corollary 3.6, Corollary 3.8]{Kuo}, continuation maps satisfy the following compatibility property with limits and colimits on the time direction:

\begin{lem} \label{lem: continuation-map-lim-colim}
In the setting of Proposition \ref{prop: continuation-map}, we have the following equivalence:
\[\underset{r \to t^-}{\colim} \, \cF_r \xrightarrow{\sim} \cF_t \ \text{and} \ \cF_t \xrightarrow{\sim} \lim_{r \to t^+} \cF_r. \]
\end{lem}

One consequence of the isotopy techniques is that we can often assume generic positions using the following proposition:

\begin{prop}[{\cite[Lemma 2.10]{ZhouIsotopy}, \cite[Proposition 3.18]{Kuo}}] \label{prop:perturbation}
Let $\cF, \cG \in \Shv(X)$ be sheaves such that $\supp(\cF) \cap \supp(\cG)$ is proper. Let $\Phi$ be a positive isotopy on $S^*X$ such that the condition 
\begin{equation} \label{for: non-transversal}
\ms^\infty(\cF) \cap  \varphi_t(\ms^\infty(\cG))  = \varnothing
\end{equation}
is satisfied for $t > 0$. Then the continuation morphism induces an isomorphism
\[\hom(\cF, \cG) \xrightarrow{\sim} \hom(\cF, \cK(\Phi)|_t \circ \cG). \]
As $\cK(\Phi)$ has $\cK(\Phi^{-1})$ as a convolution inverse, we equivalently have
\[\hom(\cF, \cG) \xrightarrow{\sim} \hom(\cK(\Phi^{-1})|_t \circ \cF,  \cG). \]
\end{prop}

\begin{notation} \label{not: wrapping-shorthand} 
When performing computation, we sometimes 
write $\cG^\Phi$ instead of $\cK(\Phi)|_t \circ \cG$. Similarly, since the inverse of a positive wrapping is a negative one, we similarly use the above notation.
\end{notation}

\subsection{Wrapped sheaves}\label{sctn_wrapping}

We now give a description of the (positive) wrapping category, which is also systematically studied in \cite{Kuo}. The point is this category organizes homotopy classes of compactly supported isotopies with fixed ends, which we then quantize its objects to wrapping kernel (thus functors) and naturally define the category of wrapped sheaves. 

\begin{defn}[Category of Positive Wrappings]\label{def:wrapping-cat} For an open set $\Omega^\infty \subseteq S^*X$, the (positive) wrapping category $W(\Omega^\infty)$ is, in principle, defined as follows: 

\begin{itemize}
    \item[(1)] {Objects:} Pairs $(\varphi, [\Phi])$ where $\varphi$ is a compactly supported contactomorphism and $[\Phi]$ is a homotopy class of compactly supported isotopies, defined on a closed interval $I$, having $\varphi$ as its endpoint, in the sense of Definition \ref{positivewrapping}.
    
    \item[(2)] {1-morphisms:} A morphism $\Psi : [\Phi_0] \to [\Phi_1]$ is a positive isotopy $\Psi: S^*X \times I \to S^*X$ such that $[\Phi_1] = [\Psi \# \Phi_0]$. Composition is given by concatenation. 
    
    \item[(3)] {2-morphisms:} For $\Psi_0, \Psi_1 : [\Phi_0] \to [\Phi_1]$, a 2-morphism is a $J$-parametric family of isotopies $\Theta : S^*X \times I \times J \to S^*X$ constant near endpoints such that $\Theta(-, t, s) = \Psi_s(-, t)$ is positive for any $s \in J$ and $\Theta(-, i, s) = \Phi_i(-)$ for $i = 0, 1$.
    
    \item[(4)] {Higher morphisms:} An $n$-morphism is a $J^{n-1}$-parametric family of isotopies $\Theta : S^*X \times I \times J^{n-1} \to S^*X$ constant near endpoints.
\end{itemize}
\end{defn}

\begin{defn}[Wrapping kernels and wrapping functors]\label{wrappingkernel}
The sheaf quantization produces a functor
\[
w:W(\Omega^\infty)\longrightarrow \Shv(M),\qquad \Phi\mapsto \cK(\Phi)|_1,
\]
sending positive isotopies to kernels whose convolution induces continuation maps on $\Shv(M)$.
For an open set $\Omega \subset S^*X$, we define the \emph{wrapping functors}
\[
\wrap^+(\Omega^\infty),\ \wrap^-(\Omega^\infty):\Shv(M)\longrightarrow \Shv(M)
\]
by
\[
\wrap^+(\Omega^\infty)(F)
:= \operatorname*{colim}_{\Phi\in W(\Omega^\infty)}\bigl(\cK(\Phi)|_1\circ F\bigr),\qquad
\wrap^-(\Omega^\infty)(F)
:= \operatorname*{lim}_{\Phi\in W(\Omega^\infty)}\bigl(\cK(\Phi)|_1\circ F\bigr).
\]
For a closed subset $\Lambda\subseteq S^*X$, we write
\[
\wrap_\Lambda^+ := \wrap^+(S^*X\setminus \Lambda),\qquad
\wrap_\Lambda^- := \wrap^-(S^*X\setminus \Lambda).
\]
\end{defn}
Geometrically, $\wrap^+(\Omega^\infty)$ pushes $F$ forward by increasingly positive wrappings in $\Omega^\infty$ and takes a colimit over all such pushes; $\wrap^-(\Omega^\infty)$ does the dual limit construction.

\begin{thm}[Adjointness of wrapping functors, {\cite[Theorem 1.2]{Kuo}}]\label{thm:wrap-adjoint}
Let $\Lambda^\infty \subseteq S^*X$ be closed, $\Lambda = X \cup (\Lambda^\infty \times \real_+)$, and $\iota_{\Lambda*}:\Shv_{\Lambda}(X)\hookrightarrow \Shv(X)$ be the inclusion of sheaves. Then there is an adjunction
\[
\wrap_\Lambda^+\dashv \iota_{\Lambda*} \dashv \wrap_\Lambda^-,
\]
i.e.~$\wrap_\Lambda^+$ (resp.~$\wrap_\Lambda^-$) is the left (resp.~right) adjoint $\iota^*_\Lambda$ (resp. $\iota^!_\Lambda$) of the inclusion $\iota_{\Lambda*}$. 
\end{thm}

We now focus on giving the definition of category of wrapped sheaves in \cite{Kuo}. 
Fix a closed conic $\Lambda^\infty \subseteq S^*X$. Ideally, one first restricts to sheaves with perfect stalks whose microsupport lies in the complement of $\Lambda^\infty$ and is subanalytic singular isotropic, 
and then formally inverts all continuation maps coming from positive isotopies in $S^*X\setminus\Lambda^\infty$. 
We now spell out the details:

Let $M$ be a real analytic manifold and $\Lambda^\infty \subseteq S^*X$ a closed subset. We first define $\widetilde{\mathfrak{wsh}}_{\Lambda}(M)$ to be the small subcategory of $\Shv(M)$ generated under finite colimits and retracts by sheaves of the form $F^{\Phi}:=\cK(\Phi)|_1\circ F$, where:
\begin{enumerate}
    \item $F$ is a constructible sheaf with perfect stalks %
    such that $\ms^{\infty}(F) \cap \Lambda^\infty = \varnothing$.
    \item $\Phi$ is a contact isotopy compactly supported away from $\Lambda$.
\end{enumerate}

Recall from Definition \ref{positivewrapping} and \ref{wrappingkernel} that now every positive isotopy $\Psi$ in $W(\Omega)$ induces a unique continuation map $c(\Psi, F) : F^{\Phi_0} \to F^{\Phi_1}$. To encode the effect of wrapping, we define the localizing subcategory:
\[
\mathcal{C}_{\Lambda}(M) := \langle \text{cof}(c(\Psi, F)) \mid \Psi \in \text{Mor}(W(S^*X \setminus \Lambda)), F \in \widetilde{\mathfrak{wsh}}_{\Lambda}(M) \rangle
\]

\begin{defn}[wrapped sheaves {\cite[Definition 4.1]{Kuo}}] The category of wrapped sheaves $\mathfrak{wsh}_{\Lambda}(M)$ is then defined as the quotient category:
\[
\mathfrak{wsh}_{\Lambda}(X) := \widetilde{\mathfrak{wsh}}_{\Lambda}(X) / \mathcal{C}_{\Lambda}(X) \cong \text{cof}\left( \mathcal{C}_{\Lambda}(X) \hookrightarrow \widetilde{\mathfrak{wsh}}_{\Lambda}(X) \right)
\]
where the cofiber is taken in stable categories. %
\end{defn}

In this localized category, morphisms can be computed as a colimit \cite[Proposition 4.3]{Kuo}:
\[
\text{Hom}_{\mathfrak{wsh}}(G, F) = \text{colim}_{F \to F^w} \text{Hom}(G, F^w)
\]
where the colimit is taken over all positive wrappings.


To conclude the section, we present the main theorem in \cite{Kuo}, where the second author showed that any compact object in $\Shv_\Lambda^c(M)$ can be represented by the wrapping of a constructible sheaf with perfect stalk avoiding the stop. 
\begin{thm}[{\cite[Theorem 1.3]{Kuo}}]\label{wrap_compare}
    Let $M$ be closed and $\Lambda \subseteq S^*X$ be a subanalytic singular Legendrian. The positive wrapping functor in Theorem \ref{thm:wrap-adjoint} satisfies $\iota_\Lambda^* = \wrap_\Lambda^+: \mathfrak{wsh}_\Lambda(X) \to {\Shv}_\Lambda^c(X)$ and induces an equivalence.
\end{thm}


\section{Compact-object-preserving kernels}
We begin the discussion of the main topic of this paper in this section. As before, we assume $X$, $Y$, etc., to be real analytic manifolds, and $\Lambda \subseteq T^* X$, $\Sigma \subseteq T^* Y$, etc. to be subanalytic.


\begin{defn} \label{defn: convolution}
Let $\Lambda \subseteq T^* X \times T^* Y$ and $\Sigma \subseteq T^* Y \times T^* Z$ be subsets. We define their composition to be
\[\Sigma \circ \Lambda \coloneqq p_{T^*X \times T^* Z} \left(  p_{T^* Y^a \times T^* Z}^{-1}(\Sigma) \cap p_{T^* X \times T^* Y}^{-1} (\Lambda) \right), \]
where $p_{T^* Y^a \times T^* Z}$ is given by $p_{T^* Y^a \times T^* Z}(x, \xi, y, \eta, z, \zeta) = (y, - \eta, z, \zeta)$. In other words, the composition $\Sigma \circ \Lambda$ is the projection of the set
\begin{align*}
    &p_{T^* X \times T^* Z} \left( \left( \Lambda \times 0_Z + 0_X \times \Sigma\right) \cap T^* X \times 0_Y \times T^* Z   \right) \\
    =&\{(x,\xi, y, \eta, z, \zeta) \mid (x, \xi, y, \eta) \in \Lambda, (y, -\eta, z, \zeta) \in \Sigma \}.
\end{align*}
\end{defn}

The main reason for introducing composition is that, when a certain transversality and properness condition is satisfied, the microsupport of convolution and its right adjoint can be measured via the composition. These results are quite standard and go back to \cite[Equation (1.6)]{GKS}, which follows from a combination of \cite[Propositions 5.4.4, 5.4.5 \& 5.4.14]{KS13}. However, 
note that the sign issue for the right adjoint $\Hom^\circ$ is a bit more confusing.

\begin{lem}[{\cite[Equation (1.6)]{GKS}}] \label{lem: convolution-estimation}
Let $L_{xy} \in \Shv(X \times Y)$ and $L_{yz} \in \Shv(Y \times Z)$. Assume 
\begin{enumerate}
\item $\pi_{xz}: X \times Y \times Z \rightarrow X \times Z$ is proper on $X \times \supp(L_{yz}) \cap \supp(L_{xy}) \times Z$. 
\item $\left( 0_X \times \ms(L_{yz}) \right) \cap \left(\ms(L_{xy}) \times 0_Z \right) \subseteq 0_{X \times Y \times Z}$.
\end{enumerate}
Then,  the sheaf $L_{yz} \circ L_{xy}  \in \Shv(X \times Z)$ satisfies 
\[ \ms\left( L_{yz} \circ L_{xy} \right) \subseteq  \ms(L_{yz}) \circ \ms(L_{xy}). \]

Respectively, let $L_{xz} \in \Shv(X \times Z)$ and $L_{yz} \in \Shv(Y \times Z)$.\footnote{Note we use the notation $x, y, z$ instead of $1, 2, 3$, since the former is less tied to a fixed ordering. That is, this convention allows us to implicitly swap $L_{xy} \subseteq T^* X \times T^* Y$ for $L_{yx} \subseteq  T^*Y \times T^* X$. } Assume 
\begin{enumerate}
\item $\pi_{xy}: X \times Y \times Z \rightarrow X \times Y$ is proper on $X \times \supp(L_{yz}) \cap Y \times \supp(L_{xz})$. 
\item $\left( 0_X \times \ms(L_{yz}) \right) \cap \left( 0_Y \times \ms(L_{xz}) \right) \subseteq 0_{X \times Y \times Z}$.
\end{enumerate}
Then, the sheaf $\Hom^\circ(L_{yz}, L_{xz}) \in \Shv(X \times Y)$ satisfies  
\[ \ms\left( \Hom^\circ(L_{yz}, L_{xz}) \right) \subseteq  (-\ms(L_{zy})) \circ \ms(L_{xz}) \subseteq T^* X \times T^* Y .\]
\end{lem}


While assumptions are needed for Lemma \ref{lem: convolution-estimation}, they can be relaxed in the following standard microsupport estimate for sheaves on a product space, which we will need later:

\begin{lem}[{\cite[Proposition 4.5]{KL25}}]\label{lem: pushtoLambda} 
    Let $\cK \in \Shv(X \times Y)$, $\cF \in \Shv(X)$ and $\cG \in \Shv(Y)$. If $\ms(\cK) \subseteq -\Lambda \times \Sigma$, $\ms(\cF) \subseteq \Lambda$ and $\ms(\cG) \subseteq \Sigma$ all of which are closed conic subsets, then 
    \[\ms(\cK \circ \cF) \subseteq \Sigma, \quad \ms(\Hom^\circ(\cK,\cG) ) \subseteq \Lambda.\]
\end{lem}

Recall, from Theorem \ref{thm: adjoint}, we use the notation $\iota_\Lambda^*$ (resp. $\iota_\Lambda^!$) to denote the left (resp. right) adjoint of the inclusion $\iota_{\Lambda*}: \Shv_\Lambda(X) \rightarrow \Shv(X)$. The above Lemma \ref{lem: pushtoLambda} implies the following compatibility between convolution and the left adjoints:

\begin{lem}[{\cite[Lemma 4.6]{KL25}}]\label{lem: niceformwrapping}
    Let $X$ and $Y$ be compact, $\Lambda \subseteq T^*X$ and $\Sigma \subseteq T^*Y$ be closed conic subsets, $\mathcal{F}\in\Shv (X)$ and $\cK\in\Shv(X\times Y) $. The following canonical maps
    \[\iota_{\Sigma}^*(\cK\circ \iota_\Lambda^*\mathcal{F}) \xrightarrow{\sim} (\iota_{-\Lambda \times \Sigma}^* \cK) \circ \iota_{\Lambda}^*\mathcal{F} \xleftarrow{\sim}  (\iota_{-\Lambda \times \Sigma}^* \cK) \circ \mathcal{F}\]
  are both isomorphisms. Passing to right adjoints, when viewing them as functors acting on $\cF$, we thus have, for $\cG \in \Shv(Y)$, 
  \[ \iota_\Lambda^! \Hom^\circ(\cK, \iota_\Sigma^! (\cG)) \xleftarrow{\sim} \iota_\Lambda^! \Hom(\iota_{-\Lambda \times \Sigma}^*(\cK), \cG) \xrightarrow{\sim} \Hom^\circ(\iota_{-\Lambda \times \Sigma}^*(\cK), \cG).\]
\end{lem}

One often finds in nature kernels $\cK \in \Shv(X \times Y)$ which are not microsupported in $-\Lambda \times \Sigma$; Example \ref{ex: cop-not-c} is such an instance. Even when the condition fails, the following argument says that one can perturb the microsupport without changing the wrapped sheaf.

Since we made the assumption that $\ms(\cK)$ is subanalytic, by \cite[Proposition 8.3.10]{KS13}, $\ms(\cK)$ is a conic Lagrangian subset. By \cite[Proposition 8.3.12]{KS13}, one can always displace subanalytic conic Lagrangians from each other (see also the explanation in \cite[Section 4.1]{KL22} \cite[Section 3.3]{LiNadlerShende}). Consequently, we have:

\begin{lem} \label{lem: perturbing-kernels}
For any $\cK \in \Shv(X \times Y)$ such that $\ms(\cK)$ is a subanalytic Lagrangian, there exists a negative wrapping $\Psi$ on $S^*(X \times Y)$ such that, for small $t > 0$,
\[ \psi_t(\dot{\ms}(\cK)) \cap \left(-\Lambda \times \Sigma\right) = \varnothing.\]
Then, we have $\iota^*_{-\Lambda \times \Sigma}(\cK) = \iota^*_{-\Lambda \times \Sigma}( \cK(\Psi) |_t \circ \cK)$.
\end{lem}
\begin{proof}
By the deformation lemma, Proposition \ref{prop:perturbation}, for any $\mathcal H \in \Shv_{-\Lambda \times \Sigma}(X \times Y)$, we have
\[ \hom(\iota^*_{-\Lambda \times \Sigma}(\cK),\mathcal H) \xrightarrow{\sim} \hom(\cK,\mathcal H) \xrightarrow{\sim} \hom( \cK^\psi, \mathcal H)  \xleftarrow{\sim} \hom(\iota^*_{-\Lambda \times \Sigma}(\cK^\psi), \mathcal H)  \]
as desired (where $\cK^\psi$ is the negative wrapping of $\cK$ by $\psi_t$ for a small $t > 0$ by Notation \ref{not: wrapping-shorthand}).  
\end{proof}

We thus need to consider the included kernel $\iota_{-\Lambda \times \Sigma}^*(\cK)$ in order to use the classification result, Theorem \ref{thm:kernel}. Then Lemmas \ref{lem: convolution-estimation} and \ref{lem: niceformwrapping} tell us that we can factorize $\iota_{-\Lambda \times \Sigma}^*(\cK) \circ (-): \Shv_{\Lambda}(X) \rightarrow \Shv_\Sigma(Y)$ to
\[
\begin{tikzcd}[column sep=0.9in]
  \Shv_\Lambda(X)
    \arrow[r, "\cK \circ (-)"] 
& \Shv_{L \circ \Lambda}(Y) 
    \arrow[r, "\iota_\Sigma^* \circ \iota_{L \circ \Lambda*}"] 
& \Shv_\Sigma(Y).
\end{tikzcd}
\]
where we set $L \coloneqq \ms(\cK)$ and the second arrow is the restriction of $\iota_\Sigma^*$ along the natural inclusion $\iota_{L \circ \Lambda*}:\Shv_{L \circ \Lambda}(Y) \hookrightarrow \Shv(Y)$. The goal is then to find geometric conditions so that each arrow is compact-object-preserving. Alternatively, one considers the functor
$\Hom^\circ(\iota_{-\Lambda \times \Sigma}^*(\cK), -) : \Shv_\Sigma(Y) \rightarrow \Shv_\Lambda(X)$, 
the right adjoint of $\iota_{-\Lambda \times \Sigma}^*(\cK) \circ (-)$,  and it factors as 
\begin{equation} \label{for: factorization}
\begin{tikzcd}[column sep=0.75in]
  \Shv_\Sigma(Y)
    \arrow[r, "\Hom^\circ(\cK\text{,}-)"] 
& \Shv_{ (-L^v) \circ \Sigma}(X) 
    \arrow[r, "\iota_\Lambda^! \circ \iota_{(-L^v) \circ \Sigma*}"] 
& \Shv_\Lambda(X).
\end{tikzcd}
\end{equation}
Here, we use $L^v$ to denote the set in $T^* Y \times T^* X$ that corresponds to $L$ in $T^* X \times T^* Y$. By Lemma \ref{lem:compact}, it is sufficient to find conditions so that the two arrows are colimit-preserving. 

\subsection{Colimit preserving internal Hom} \label{ssec: internal-Hom}
We begin with considering the functor $\Hom^\circ(\cK,-)$. In this subsection, we will show the following sufficient condition for  $\Hom^\circ(\cK,-)$ to preserve colimits:

\begin{thm}\label{thm:colim-preserve-hom}
Let $X, Y$ be real analytic manifolds and $\Lambda \subseteq T^*X, \Sigma \subseteq T^*Y$ be subanalytic conic Lagrangians. Let $\cK \in \Shv^b_{\real-c}(X \times Y)$ be a sheaf kernel such that
\[\ms(\cK) \cap \Lambda \times \Sigma \subseteq 0_{X \times Y} \quad \text{and} \quad -\ms(\cK)^v \circ \Sigma \subseteq \Lambda.\]
Then the corresponding right adjoint integral transform preserves colimits:
\[\Hom^\circ(\cK,-): \Shv_\Sigma(Y) \to \Shv_{\Lambda}(X).\]
\end{thm}

For this purpose, we recall that for any $X$ in a closed symmetric monoidal category $\cC$, one can always take the naive dual $\inhom(X,1_\cC)$, which always comes with an evaluation map
\begin{equation} \label{for: evaluation}
\inhom(X,1_\cC) \otimes Y \rightarrow \inhom(X,Y),
\end{equation}
for any $Y \in \cC$. Furthermore, $X$ is dualizable if and only if the map in Formula \eqref{for: evaluation} is an isomorphism for all $Y$.

Note that the virtue of this identification, when $\cC$ is a presentable symmetric monoidal category, is that it implies that $\inhom(X,Y)$ is both limit-preserving and colimit-preserving in $Y$. For our coefficients, dualizability coincides with being compact, so we require our sheaves to have perfect stalks. However, there are two extra conditions needed: constructibility and transversality.

\begin{prop}[{\cite[Proposition 5.4.14]{KS13}}]\label{prop:tensor-hom}
    Let $\cF, \cG \in \Shv(X)$ be sheaves such that $\ms^\infty(\cF) \cap \ms^\infty(\cG) = \varnothing$ and suppose $\cF$ is constructible with perfect stalks. Then there is a natural isomorphism
    \[\Hom(\cF, 1_X) \otimes \cG \xrightarrow{\sim} \Hom(\cF, \cG).\]
\end{prop}

The following result is standard and a related proposition can be found in \cite[Proposition 1.10]{GKS}. For the reader's convenience as well as matching notations, we still give the proof:

\begin{prop} \label{prop: right-adjoint-kernel}
Let $X$ and $Y$ be manifolds and let $\cK \in \Shv(X \times Y)$ be a sheaf kernel. Assume
\begin{enumerate}
\item The projection to $X$, $p_X: X \times Y \rightarrow X$, is proper on $\supp(\cK)$,
\item $\cK$ is constructible with perfect stalks, and
\item $\ms(\cK) \cap 0_X \times \Sigma \subseteq 0_{X \times Y}$.
\end{enumerate}
Then, $\Hom^\circ(\cK, -): \Shv_\Sigma(Y) \rightarrow \Shv_{ (-L^v) \circ \Sigma}(X) $ is equivalent to $\cK^R \circ(-)$
where \[\cK^R \coloneqq \Hom(v^* \cK, 1_Y \boxtimes \omega_X)\]
and $v: X \times Y \cong Y \times X$ is the coordinate swapping map $v(x,y) = (y,x)$. 
\end{prop}

\begin{proof}
Let $G \in \Shv(Y)$. Our goal is to turn ${p_{X*}} \Hom(\cK, p_Y^! G)$ into a convolution. By  \cite[Proposition 5.4.5]{KS13}, $\ms(p_Y^! G)$ is microsupported in $0_X \times \ms(G)$ and 
$p_Y^! G = p_Y^* G \otimes p_Y^! 1_Y$. 
Assumptions (2) and (3) then imply that Proposition \ref{prop:tensor-hom} is applicable and 
\[\Hom(\cK, p_Y^! G) = \Hom(\cK, 1_{X \times Y}) \otimes p_Y^! G. \]
To get the correct form for convolution, by the proof of \cite[Proposition 3.3.2]{KS13}, we have that $p_Y^! 1_Y = p_X^* \omega_X = \omega_X \boxtimes 1_Y$. Lastly, Assumption (1) implies that 
\[{p_{X*}} \Hom(\cK, p_Y^! G) 
= {p_{X!}}  \left(\Hom(v^* \cK,  1_Y \boxtimes \omega_X ) \otimes p_Y^* G \right),\]
as desired. This completes the proof.
\end{proof}

\begin{proof}[Proof of Theorem \ref{thm:colim-preserve-hom}]
It follows from Proposition \ref{prop: right-adjoint-kernel} that under the given assumptions, we have $\Hom^\circ(\cK, -) \cong \cK^R \circ (-)$. We know the convolution is defined by left adjoints, and thus preserves colimits.
\end{proof}

A closely related result also appears in \cite[Section 8.4]{LiNadlerShende}. In the case of sheaves on manifolds, we have the following theorem. 

\begin{prop}\label{thm:wrap-one-factor} 
Let $X, Y$ be   real analytic manifolds and $\Lambda \subseteq T^*X, \Sigma \subseteq T^*Y$ be subanalytic conic Lagrangian subsets. Let $\cK \in \Shv_{\mathbb{R}-c}^b(X \times Y)$ be a sheaf kernel, and suppose for
$$\ms(\cK) \cap \Lambda \times \Sigma \subseteq 0_{X \times Y} \quad \text{and} \quad -\ms(\cK)^v \circ \Sigma \subseteq \Lambda.$$ 
Then for any 
sheaf $\cF \in \Shv(X)$, 
$$\iota_{-\Lambda \times \Sigma}^*\cK \circ \iota_\Lambda^*\cF = \iota_\Sigma^*(\cK \circ \cF).$$
\end{prop}
\begin{proof}
	By Lemma \ref{lem: niceformwrapping}, we in general only have $(\iota_{-\Lambda \times \Sigma}^* \cK) \circ \iota_{\Lambda}^*\cF = \iota_{\Sigma}^*(\cK\circ \iota_\Lambda^*\cF)$ and the goal is to show that the right hand side is $\iota_\Sigma^*(\cK \circ \cF)$ under the assumption $- {\ms}^\infty(\cK)^v \circ \Sigma \subseteq \Lambda$. Indeed, let $\cG \in \Shv_\Sigma(Y)$ and we can compute that
	\[
	\hom(\iota_{\Sigma}^*(\cK\circ \iota_\Lambda^*\cF), \cG) = \hom( \cK\circ \iota_\Lambda^*\cF, \cG) = \hom( \iota_\Lambda^* \cF, \Hom^\circ(\cK, \cG)). 
	\]
	But then the assumption, by Lemma \ref{lem: pushtoLambda}, implies that 
	\[ \hom( \iota_\Lambda^* \cF, \Hom^\circ(\cK, \cG)) = \hom(  \cF, \Hom^\circ(\cK, \cG)) = \hom(\iota_\Sigma^*(\cK \circ \cF), \cG). \]
    This completes the proof.
\end{proof}


\begin{proof}[Alternative proof of Theorem \ref{thm:colim-preserve-hom}]
Let $\cF \in \Shv^c_{\Lambda}(X)$. Then there exists a constructible sheaf with perfect stalk and compact support $\cF_0 \in \Shv^b_{\real-c}(X)_0$ such that $\cF = \iota_\Lambda^* \cF_0$.
As $\cK \circ \cF_0 = {\pi_Y}_! (\cK \otimes \pi_X^* \F)$ and we assume the projection to $Y$ is proper, the stalk
\[(\cK \circ \cF_0)_y = \Gamma_c(X \times \{y\}; \cK|_{X \times \{y\}} \otimes F) = \Gamma(X \times \{y\}; \cK|_{X \times \{y\}} \otimes F), \]
being the global sections of a constructible sheaf with perfect stalks, is perfect. By Theorem \ref{wrap_compare}, this also implies that $\iota_{-\Lambda \times \Sigma}^*(\cK) \circ  \cF = \iota_\Sigma^*(\cK \circ \cF) \in \Shv^c_\Sigma(Y)$.
\end{proof}

\subsection{Colimit preserving negative wrapping} \label{ssec: squeezability}
Thus it boils down to give a geometric condition that makes the functor
\begin{equation} \label{for: iota-shriek} 
\iota_\Lambda^! |_{(-L^v) \circ \Sigma}: \Shv_{ (-L^v) \circ \Sigma}(X) \rightarrow \Shv_\Lambda(X),
\end{equation}
which is a priori only a right adjoint, colimit-preserving. Note that since we rely on Lemma \ref{lem: perturbing-kernels} for condition (2) in Proposition \ref{prop: right-adjoint-kernel} to be satisfied, we need the condition to be open. The following is one of those:

\begin{defn} \label{def: squeezable}
Let $V$ be a co-oriented contact manifold. Denote the contact isotopy on $V$ induced by the Reeb flow by $T_t$. Let $\Lambda \subseteq V$ be a compact Legendrian, and set $\Lambda_\epsilon \coloneqq T_\epsilon(\Lambda)$. We say that a compact Legendrian $\Gamma \subseteq V$ with $\Gamma \cap \Lambda = \varnothing$ is negatively squeezable to $\Lambda$ with bounded energy if there exists a negative isotopy $\Psi: V \times [0,1) \rightarrow V$ such that 
\begin{enumerate}
\item $\psi_t |_\Lambda = \id_\Lambda$ for all $t \in [0, 1)$.
\item For any neighborhood $\cU$ of $\Lambda_\epsilon$, there exists $T \in [0,1)$ such that $\psi_t(\Gamma) \subseteq \cU$ for all $t \in [T, 1)$.
\item There exists $C > 0$ such that its Hamiltonian $H_\Psi$ satisfies $|H_\Psi|_{\psi_t(\Gamma) \times \{t\}}| < C$ for all $t \geq 0$.
\end{enumerate}
\end{defn}

\begin{ex} \label{ex: Liouville-hypersurface}
Let $\Lambda \subseteq S^*X$ be the skeleton of a Liouville hypersurface $F \subseteq S^*X$, and let $\Gamma \subseteq S^*X$. Suppose there is a contact isotopy $\Phi$, supported away from $\Lambda$, that sends $\Gamma$ into a tubular neighborhood $F \times (\epsilon, 2\epsilon)$ at time $1$. Then $\Gamma$ is squeezable into $\Lambda$ with bounded energy.

In fact, in this case, one can consider the contact form on $F \times (\epsilon, 2\epsilon)$ given by $dz + \lambda_F$. A contact lift of the Liouville flow $\Phi^Z$ is generated by the Hamiltonian $H(p, z) = \epsilon - z$. We can then define a contact isotopy $\Psi$ supported away from $\Lambda$, such that $\psi_t = \phi_t$ for $t \in [0, 1]$ and
\[\psi_t = \phi^Z_{t-1} \circ \phi_1, \quad t \in [1,\infty).\]
Then $\Psi$ is generated by the Hamiltonian $H_\psi(p, z) = \epsilon - z \circ \phi_1(p, z)$ for $t \geq 1$, where $z \in (\epsilon, 2\epsilon)$, and it satisfies all the conditions above.

For example, to see the squeezing, we first simplify the notation and assume that $\Gamma$ is contained in $F \times (\epsilon, 2\epsilon)$ from the beginning. Following \eqref{for: contact-flow}, the vector field generated by the Hamiltonian $H = \epsilon - z$ is $V_H = -Z_F + (\epsilon - z)\frac{\partial}{\partial z}$, where $Z_F$ is the Liouville vector field on $F$. Integrating, the $z$-component of the flow is given by $z(t) = (z(0) - \epsilon)e^{-t} + \epsilon$, so $\lim_{t \to \infty} z(t) = \epsilon$, as desired.

To see the bounded energy condition, consider the reparametrization 
\[ \tau = \frac{t}{t+1}, \quad \text{so} \quad t = \frac{\tau}{1-\tau} \quad \text{and} \quad \frac{dt}{d\tau} = \frac{1}{(1-\tau)^{2}}, \]
so that the time direction is now parametrized by $[0,1)$. We check that, under this reparametrization, the Hamiltonian is bounded. Note that, while the reparametrization also changes the speed of the Liouville flow on $F$, this effect is invisible when pairing with $dz + \lambda_F$. Thus we only need to consider the change in the $z$-direction. Since $H_{\psi}(t) = \alpha_{\psi_t}(\partial_t \psi_t)$, by the chain rule,
\[\tilde{H}_{\psi}(\tau) = \alpha_{\psi_{t(\tau)}}(\partial_\tau \psi_{t(\tau)}) = \alpha_{\psi_{t(\tau)}}(\partial_t \psi_{t}) \cdot \frac{d t}{d \tau} = H_{\psi}(t) \cdot \frac{d t}{d \tau},\]
so, under the $\tau$-parametrization, $\tilde{H}_\psi = (\epsilon - z)(1-\tau)^{-2}$. Since $z(t) - \epsilon = (z(0) - \epsilon)e^{-t}$, we compute (by L'H\^{o}pital's rule) that
\[|\tilde{H}_\psi|_{\psi_\tau(\Gamma)}| \leq C e^{-\tau/(1-\tau)}(1-\tau)^{-2} \to 0, \quad \tau \to 1.\]
Hence $\tilde{H}_\psi$ is bounded on $[0,1)$.
\end{ex}

\begin{ex} \label{ex: swappable}
If one defines similarly in Definition \ref{def: squeezable} the notion of positively squeezable, then we recover the notion of swappability first studied by Sylvan in \cite[Definition 1.1]{Sylvan-Orlov} in the wrapped Fukaya category framework, and later by the second and third authors in \cite[Definition 5.19]{KL22}. More precisely, a compact Legendrian $\Lambda \subseteq S^* X$ is swappable, by definition, if and only if $\Lambda_\epsilon$ is positively squeezable to $\Lambda$ and $\Lambda_{-\epsilon}$ is negatively squeezable to $\Lambda$.
\end{ex}

\begin{prop} \label{prop: squeezable-iota-shriek} 
Let $X$ be a manifold and $\Lambda, \Gamma \subseteq S^* X$ be subanalytic Legendrians. If $\Gamma$ is negatively squeezable to $\Lambda$ with bounded energy, then the functor
\[\iota_\Lambda^!: \Shv_{\Gamma}(X) \to \Shv_\Lambda(X)\]
is colimit-preserving; respectively, if $\Gamma$ is positively squeezable to $\Lambda$, 
the functor
\[\iota_\Lambda^*: \Shv_{\Gamma}(X) \to \Shv_\Lambda(X)\]
is limit-preserving.
\end{prop}

We will need some technical results for this proposition. First, by \cite[Lemma 5.7]{Kuo}, for $i: X \times \{\infty\} \hookrightarrow X \times [0, +\infty]$ and $j: X \times [0,+\infty) \hookrightarrow X \times [0, +\infty]$, we have, for $\cF \in \Shv_\Gamma(X)$,
\[\iota_\Lambda^!(\cF) = \wrap_\Lambda^-(\cF) = i^*j_*(\cK(\Psi) \circ \cF).\]
That is, as $\Psi$ pushes $\Gamma$ into $\Lambda$, instead of running through all negative wrappings compactly supported away from $\Lambda$, we can consider $\lim_{t \to \infty} (\cK(\Psi)|_{t} \circ \cF)$, which can in turn be computed as the nearby cycle functor (at $\infty$) of the entire ``family of sheaves'' $\cK(\Psi) \circ \cF$. Here, the nearby cycle functor in our setting has been studied extensively in \cite{NS20} and plays a key role in the development of microsheaves on Weinstein manifolds. 

The key claim here is that one can rewrite $i^*j_*(\cK(\Psi) \circ \cF)$, under the bounded energy assumption, as a colimit followed by the action of the time-$1$ isotopy. See \eqref{for: main-prop-end-formula}. Naively speaking, we have
\[ i^*j_*(\cK(\Psi) \circ \cF) = i^*j_* \left( \cK(T_{-Ct}) \circ|_t \cK(T_{Ct} \circ|_t \Psi)|_t \circ \cF\right),\]
so we would be done if taking nearby cycles commuted with convolution. In the special case of convolution with a GKS kernel, we can use a result of the third author, concerning the family version of the nearby cycle, to see that this is indeed the case; this is the main argument of the proof. For this purpose, we first recall a special case of the definition of relative microsupport from \cite[Section 3.2]{NS20}:

Let $X$ and $B$ be manifolds, where we view $B$ as a smooth family of parameters. The projection $\pi: X \times B \rightarrow B$ induces a vector bundle map (over $X \times B$), $\pi^* T^*B \rightarrow T^* X$. As we are in the split situation, we can trivialize it and complete it to a short exact sequence
\begin{equation} \label{for: capital-pi}
0 \rightarrow 0_X \times T^* B \rightarrow T^* (X \times B) \xrightarrow{\Pi} T^* X \times 0_B \rightarrow 0.    
\end{equation}

\begin{defn}[{\cite[Definition 3.10]{NS20}}]
For $\mathcal H \in \Shv(X \times B)$, the relative microsupport of $H$ is the set 
\[ \ms_\pi(\mathcal H) \coloneqq \overline{\Pi \left( \ms(\mathcal H)\right)} \subseteq T^* X \times 0_B. \]
\end{defn}

\begin{prop}[{\cite[Proposition 3.3]{LiCobordism}}] \label{prop: family-base-change}
Let $X$ be a manifold, and consider $a, b, \epsilon > 0$ \footnote{The choice of these numbers is not essential. We make this choice only to avoid having a manifold with boundary when imposing microsupport conditions. The same applies to the choice of $(-1,2)$ in \eqref{for: perturb-by-Reeb}.} and  the commutative diagram
\[
\begin{tikzcd}
	X \times (-\epsilon,+\infty) \times \{0\} && X \times (-\epsilon,+\infty) \times (-a,b)  \\
	X \times (-\epsilon, +\infty] \times \{0\} && X \times (-\epsilon, +\infty] \times (-a,b)
	\arrow["i", from=1-1, to=1-3]
	\arrow["\bar{i}", from=2-1, to=2-3]
	\arrow["j", from=1-1, to=2-1]
	\arrow["\bar{j}", from=1-3, to=2-3]
\end{tikzcd}
\]
so, for $\mathcal H \in \Shv(X \times (-\epsilon,+\infty) \times (-a,b))$, there exists a canonical base change map
\begin{equation} \label{for: family-base-change}
   \bar{i}^*\bar{j}_* \mathcal H \rightarrow  j_* i^* \mathcal H.
\end{equation}
Then, the map \eqref{for: family-base-change} is an equivalence if $\mathcal H$ satisfies the following conditions:\footnote{In the original result, the first non-characteristic condition does not include fibers at $0 \in (-a, b)$, and instead the author used a non-characteristic condition of $i^\sharp \ms(\mathcal H)$ where $i^\sharp$ is an enlarged upper bound used to estimate the effect of $i^*$ in the characteristic situation. However, our case is in fact noncharacteristic, and the usual microsupport estimation suffices \cite[Remark 6.2.6]{KS13}.} 
\begin{enumerate}
    \item $\ms(\mathcal H) \cap \left(0_X \times T^* (-\epsilon,+\infty) \times T^* (-a,b) \right) \subseteq 0_{X \times (-\epsilon,\infty) \times \left(-a,b \right)}$.
    \item $\ms_{\pi}(\mathcal H) \cap (T^* X \times \{0,0\})$ is a subanalytic conic Lagrangian.
\end{enumerate}   
\end{prop}

When applying the above proposition, we think of $\mathcal H$ as an $(-a,b)$-family of time-dependent sheaves, and we may think of the effect of $j_*$ and $\bar{j}_*$ as extending the sheaf to $+\infty$ by ``taking the limit". We thus use the notation $\mathcal H_0 \coloneqq \mathcal H |_{X \times (-\epsilon,+\infty) \times \{0\}}$, $\mathcal H_0(+\infty) \coloneqq (j_* \mathcal H_0)|_{X \times \{+\infty\}}$, and $\mathcal H(+\infty)  \coloneqq (\bar{i}^* \bar{j_*} \mathcal H)|_{X \times \{+\infty\} \times (-a,b) }$. Restricting \eqref{for: family-base-change} to $X \times \{+\infty\} \times \{0\}$, we see that the limit at $+\infty$ of the $s =0$-slice can be computed by first taking the limit of the entire family and then setting $s= 0$. More precisely,

\begin{cor} \label{cor: family-base-change}
    Let $\mathcal H \in \Shv(X \times (-\epsilon,+\infty) \times (-a,b))$ be a sheaf satisfying the condition in Proposition \ref{prop: family-base-change}. Then the canonical map 
    \[\mathcal H_0(+\infty) \xleftarrow{\sim} \mathcal H(+\infty) |_{s = 0}. \]
\end{cor}

We will apply the corollary in the situation where $\mathcal H$ is a family-isotopy of a sheaf: let $\Psi:\dT^* X \times I \rightarrow \dT^* X$ be a negative isotopy squeezing $\Gamma$ into $\Lambda$. For convenience, although we stated Corollary \ref{cor: family-base-change} with the convention that time $t$ is parametrized by $\real$, we will assume here that $t \in I = (-1,1)$ instead, so that $\psi_t(\Gamma)$ is squeezed into $\Lambda$ as $t \to 1$.

Now, we define a family of isotopies
\begin{align} \label{for: perturb-by-Reeb}
    \Theta: \dT^* X \times (-1,1) \times (-1,2) &\rightarrow \dT^* X \\
    \Theta(x,\xi,t,s) &= T_{Cst} \circ \psi_t(x,\xi).\notag 
\end{align} 
where $C > 0$ is a fixed constant such that $T_{Ct} \circ \psi_t$ is a positive isotopy. Recall, as in Definition \ref{def: squeezable}, we use $T_t$ to denote the contact isotopy generated by the Reeb vector field $R_\alpha$. 

\begin{lem}\label{lem:hamiltonians-of-Theta}
Let $(V,\alpha)$ be a co-oriented contact manifold with Reeb vector field $R_\alpha$,
and let $T_c$ denote the time-$c$ flow of $R_\alpha$. Then:
\begin{enumerate}
    \item For fixed $s$, the isotopy $t \mapsto \Theta_{s,t}$ is generated by the
    contact Hamiltonian
    \[
    \alpha_{\varphi_{t,s}(x,\xi)}(\partial_t\varphi_{t,s}(x,\xi)) = Cs + H_\psi \circ T_{-Cst}.
    \]
    \item For fixed $t$, the isotopy $s \mapsto \Theta_{s,t}$ is generated by the
    constant contact Hamiltonian
    \[
    \alpha_{\varphi_{t,s}(x,\xi)}(\partial_s\varphi_{t,s}(x,\xi)) = Ct.
    \]
\end{enumerate}
\end{lem}

\begin{lem} \label{lem: F-Psi-properties}
Let $\Theta$ be defined as in Formula \eqref{for: perturb-by-Reeb}. Per Remark \ref{rmk: family-GKS}, we denote by
\[\cK(\Theta) \in \Shv(X \times X \times (-1,1) \times (-1,2))\] 
the GKS kernel of $\Theta$, and set, for $\cF \in \Shv_\Gamma(X)$,
\[\cF^\Theta = \cK(\Theta) \circ \cF \in \Shv(X \times (-1,1) \times (-1,2)).\]
Then $\cF^\Theta$ satisfies the assumptions in Proposition \ref{prop: family-base-change}. 
\end{lem}

\begin{proof}
Convolving with the GKS kernels from Theorem \ref{thm:gks} always satisfies the assumption in Lemma \ref{lem: convolution-estimation}, so we have 
\begin{equation} \label{for: F-Psi-estimation}
 \ms^\infty(\cF^\Theta) \subseteq \Lambda_\Theta \circ \ms^\infty(\cF) = \{((T_{Cst} \circ \psi_t) (x, \xi), t, -Cs - H_{\psi} \circ T_{-Cst}, s, -Ct) \mid (x, \xi) \in \ms^\infty(\cF)\}    
\end{equation}
where $(x, \xi) \in \ms^\infty(\cF)$, $t\in (-1,1)$, and $s \in (-1,2)$, so condition (1) is immediate. 
Lastly, for condition (3), we again use Formula \eqref{for: F-Psi-estimation} to see that, away from the zero section, we have
\[\ms_\pi(\cF^\Psi) = \overline{\Pi(\ms(\cF^\Psi))} \subseteq \{((T_{Cst} \circ \psi_t) (x, \xi), t, s) \mid (x, \xi) \in \ms^\infty(\cF)\} \]
so, at $(0,0)$, $\ms^\infty(\cF) \subseteq \Gamma$ and it is a subanalytic Lagrangian.
\end{proof}

\begin{proof}[Proof of Proposition \ref{prop: squeezable-iota-shriek} ] \label{proof: proof-of-main-thm}
Let $\cF \in \Shv_\Gamma(X)$. By Lemma \ref{lem: F-Psi-properties}, we can apply Corollary \ref{cor: family-base-change} to $\cF^\Psi$. As explained right below the statement of the current Proposition \ref{prop: squeezable-iota-shriek}, $\mathcal H_0(+\infty) \eqqcolon i^* j_* (\cK(\Psi) \circ \cF) = \iota_\Lambda^!(\cF)$. Thus, our goal is to show that $\mathcal H(+\infty)|_{s = 0} \coloneqq \bar{i}^*\bar{j}_* \cF$ has an expression that is colimit-preserving in $\cF$. For this purpose, we use \cite[Definition 3.17]{NS20} to estimate its microsupport:
\[\ms(\bar{i}^* \bar{j}_* \cF^\Psi) \subseteq \overline{\Pi_1(\ms(\cF^\Psi))} \cap T^* X \times \{1\} \times T^* (-1,2).\]
Here, we apply the same construction as in Formula \eqref{for: capital-pi} to $\pi_1: X \times (-1,1) \times (-1,2) \rightarrow (-1,1)$, and we note that, unlike above, as we take the nearby cycle, the closure has to be taken with $(-1,1)$ extended to $(-1,1]$.

Now, for any $(y,\eta, 1, s, \sigma) \in \overline{\Pi_1(\ms(\cF^\Psi))}$, there is a sequence $(\Theta_{(t_n,s_n)}(x_n,\xi_n),t_n, s_n, -C t_n)$ converging to it. Clearly, this means $\sigma = -C$. To see where $(y,\eta)$ is, we can equivalently apply $T_{-C}$ and find that $(T_{-C}(y,\eta), 1, s, \sigma)$ is approached by $(\psi_{t_n}(x_n,\xi_n), 1, s_n, - C t_n)$. Because $\cF \in \Shv_\Gamma(X)$, the squeezable assumption implies that $T_{-C}(y,\eta) = \lim_{n \to \infty} \psi_{t_n}(x_n,\xi_n) \in \Lambda$. Translating the points back with $T_C$, we conclude that
\[\ms(\mathcal H(+\infty))^\infty \subseteq \{(T_C(y,\eta),s, -C)| (y,\eta) \in \Lambda\}\]
or, equivalently, 
\[ \ms\left( \cK(T_{-C}) \circ|_s \mathcal H(+\infty)\right)^\infty \subseteq \Lambda \times 0_{(-1,2)}.\]
Here, we recall that the notation $\circ|_s$ means the family-convolution \cite[Equation (1.13)]{GKS}. As there is no microsupport on the second factor, by \cite[Corollary 1.6]{GKS}, this family is a constant family; that is,
\[\cK(T_{-C}) \circ|_s \mathcal H(+\infty) = \mathcal H_0(+\infty) \boxtimes 1_{(-1,2)}.\]
Restricting along $s = 1$, we obtain that
\[\mathcal H_0(+\infty) = \cK(T_{-C}) \circ \mathcal H(+\infty)|_{s = 1}.\]
To finish the proof, we note that the same condition that works for $s = 0$ also works for $s =1$, so Corollary \ref{cor: family-base-change} implies that
\[
\mathcal H(+\infty)|_{s = 1} = i^*j_*( \cK(T_{Ct} \circ|_t \Psi) \circ \cF) = \colim_{t \rightarrow 1} \cK(T_{Ct}\circ|_t \Psi) |_t \circ \cF\]
where the key fact we use here is that, unlike $\Psi$, its composition with the Reeb flow $T_{Ct}$ is chosen to be positive, so, by \cite[Lemma 5.7]{Kuo}, the nearby cycle is computed by a colimit. In summary,
\begin{equation} \label{for: main-prop-end-formula}
\iota_\Lambda^!(\cF) = \cK(T_{-C}) \circ \colim_{t \to 1} (\cK(T_{Ct} \circ|_t \Psi)|_t \circ \cF) \end{equation} is colimit-preserving.
\end{proof}

\begin{rmk}
Trading $\iota_\Lambda^!(\cF)$ for the nearby cycle of $\cK(\Psi) \circ \cF$ requires the assumption that $\Psi$ is supported away from $\Lambda$. However, when computing the latter as above, we do not need to assume that $T_{Ct} \circ|_t \Psi$ is also supported away from $\Lambda$. Therefore, $C > 0$ can be taken to be very large.
\end{rmk}

As Definition \ref{def: squeezable} is an open condition, we may assume the conditions in Proposition \ref{prop: right-adjoint-kernel} are satisfied. Combined with Proposition \ref{prop: squeezable-iota-shriek}, the factorization \eqref{for: factorization} implies the following theorem.   

\begin{thm}\label{thm:main}
 Let $X, Y$ be real analytic manifolds and $\Lambda \subseteq T^*X$, $\Sigma \subseteq T^*Y$ be conic subanalytic Lagrangian subsets. For any constructible sheaf with perfect stalks $\cK \in \Shv_{\mathbb{R}-c}^b(X \times Y)$, 
    suppose $\ms(\cK) \cap (-\Lambda \times \Sigma) \subset 0_{X \times Y}$, the projection $\supp(\cK) \subseteqq X \times Y \rightarrow X$ is proper and $-\ms(\cK)^v \circ \Sigma$ is negatively squeezable to $\Lambda$ with bounded energy (Definition \ref{def: squeezable}). Then the corresponding wrapped sheaf defines an integral transform preserving compact objects
    \[\iota_{-\Lambda \times \Sigma}^*(\cK) \circ (-): \Shv^c_{\Lambda}(X) \to \Shv^c_{\Sigma}(Y).\]
\end{thm}

\begin{rmk}
In \cite[Proposition 1.10]{KL22}, the second and third authors show that the wrap-once functor $S_\Lambda^+$, mentioned in Example \ref{ex: cop-not-c}, is an equivalence so in particular it is compact-object-preserving. Proposition \ref{prop: squeezable-iota-shriek} implies that, by Example \ref{ex: swappable}, up to the bounded energy assumption, the one-sided condition is enough to imply compact-object-preserving.
\end{rmk}

\begin{rmk} \label{rmk: looser-bounded-energy}
As one can see from the proof, the key property of $T_{Ct}$ is that $T_{Ct} \circ|_t \Psi$ is non-negative and $\lim_{t \to 1} T_{-Ct} = T_{-C}$ is a contactomorphism. Thus, we can replace the bounded energy condition (3) in Definition \ref{def: squeezable} by a looser condition. Fix $\Psi: S^* X \times [0,1) \rightarrow S^* X$; we can instead ask:
\begin{enumerate}
    \item[(3)$'$] There is a $C^\infty$ positive increasing function $\rho$ on $[0,1)$ such that $\lim_{t \to 1^-} \rho(t) = C < \infty$ and such that $T_{\rho(t)} \circ|_t \Psi$ is non-negative on $\Gamma$.
\end{enumerate}
For example, recall that in Example \ref{ex: Liouville-hypersurface}, under the coordinate $t \in [0,\infty)$, the isotopy $\psi_t$ has its $z$-component given by $z(t) = (z(0) -\epsilon)e^{-t} + \epsilon$. Take some $h$ and consider $T_{h(t)} \circ \psi_t$, so that the $z$-component is now $(z(0) -\epsilon)e^{-t} + \epsilon + h(t)$. One can then compute that the Hamiltonian of this composite isotopy is given by $H = (\epsilon + h(t) - z) + h'(t)$. Thus, any positive increasing function $h$ with $\lim_{t\to\infty} h(t) = C < \infty$ and $C > \epsilon$ gives such a $\rho$, by setting $\rho(\tau) = h(t)$ via the same change of parameter $t = \frac{\tau}{1-\tau}$ from Example \ref{ex: Liouville-hypersurface}. One concrete choice is
\[h(t) = C(1-e^{-t}).\]
\end{rmk}


\subsection{Generalization to Weinstein manifolds}

As recalled in Definition \ref{defn: microsheaf-Weinstein}, the notion of microsheaves can be extended to the setting of Weinstein manifolds. In a recent paper \cite{LiNadlerShende}, the third author, along with Nadler and Shende, has extended the notion of convolution to this setting, and thus wrappings and other constructions depending on it. This allows us to generalize Theorem \ref{thm:main} to this setting as well, which should generalize Gao's result, Theorem \ref{conj1}, through the identification with Fukaya categories \cite{GPS3}; see Lemma \ref{pre-flowon1factor} below.

Let $M, N$ be Weinstein manifolds with mostly Legendrian stops $\Lambda^\infty \subseteq \partial^\infty M, \Sigma^\infty \subseteq \partial^\infty N$. We denote by $\Lambda = \mathfrak{c}_{M,\Lambda^\infty}$ and $\Sigma = \mathfrak{c}_{N,\Sigma^\infty}$ the relative Lagrangian skeleton.

\begin{thm} \label{main-Weinstein} 
Let $M, N$ be polarizable Weinstein manifolds with mostly Legendrian stops $\Lambda^\infty \subseteq \partial^\infty M, \Sigma^\infty \subseteq \partial^\infty N$. Let $\cK \in \msh^b_{\real-c}(M^- \times N)$\footnote{Note that it does not make sense to take ``$-$'' on a general Weinstein manifold $M$, as in Theorem \ref{thm:main}. Instead, it is symplectically more natural to replace $M$ by $M^{-}$, equipped with the reversed Liouville structure.} be a microsheaf with perfect stalks and compact conic Lagrangian support. Suppose 
$\ms(\cK)^v \circ \Sigma$ is negatively squeezable to $\Lambda$ with bounded energy, the functor
\[\iota_{-\Lambda \times \Sigma}^*(\cK) \circ (-):\msh_{\Lambda}(M) \rightarrow \msh_{\Sigma}(N)\]
preserves compact objects. \footnote{The object $\cK$ lives in a category-valued sheaf, so it has a notion of support. Here we denote it by $\ms(\cK)$, since microlocally, i.e., locally on the Weinstein manifold, it is computed by the microsupport $\ms(K)$ for any representative sheaf $K$.}
\end{thm}

\begin{proof}

Recall that, by Definition \ref{defn: microsheaf-Weinstein}, microsheaves on Weinstein manifolds are realized as classical microsheaves, in the sense of Definition \ref{def:microsheaf}. In particular, $L \coloneqq \ms(\cK)$, a conic Lagrangian in $M^- \times N$, can be viewed as a subset of a cosphere bundle, and thus has an antipode, denoted $-L$. 

As in the case of sheaves, the adjunction \[\iota_{-\Lambda \times \Sigma}^*(\cK) \circ (-):\msh_{\Lambda}(M) \rightleftharpoons \msh_{\Sigma}(N): \Hom^\circ(\iota^!_{-\Lambda \times \Sigma}(\cK), -)\] implies that, by Lemma \ref{lem:compact}, showing that the functor preserves compact objects is equivalent to showing that its right adjoint, $\Hom^\circ(\iota^!_{-\Lambda \times \Sigma}(\cK), -)$, is colimit-preserving. To do this, we use the doubling functor in Proposition \ref{prop: doubling}, which realizes microsheaves as sheaves, so that we can apply the results from the previous sections.

Since $\iota^!_{-\Lambda \times \Sigma}$ commutes with the doubling functor (cf.~\cite[Corollary 5.4.12]{LiNadlerShende}), the hom-convolution with the wrapped microsheaf kernel $\iota_{-\Lambda \times \Sigma}^!(\cK)$ is realized as 
\[m_{\Lambda}\Hom^\circ( \iota_{-\Lambda_{\cup,\epsilon}^- \times \Sigma_{\cup,\epsilon}^-}^!(w_{L}^-(\cK)), w_{\Lambda}^-(-)): \msh_{\Lambda}(M) \to \msh_{\Sigma}(N).\] 
In plain words, the right adjoint of convolution of microsheaves is nothing but the procedure of realizing all the actors as sheaves via doubling, applying the usual right adjoint $\Hom^\circ$, and then microlocalizing the result.
By Lemma \ref{lem: niceformwrapping}, we can pull out $\iota^!$ from $\Hom^\circ$ and obtain the expression
\[m_{\Lambda}\iota_{\Lambda_{\cup,\epsilon}^+}^!\Hom^\circ(w_{L}^-(\cK) , w_{\Lambda}^-(-)): \msh_{\Sigma}(N) \to \msh_{\Lambda}(M),\]
so we are in a position to perform the analysis in \eqref{for: factorization}.

We first explain that the results from Section \ref{ssec: internal-Hom} go through. Here, we follow the convention of contact composition in \cite[Section 9]{LiNadlerShende}. Set $\Gamma \coloneqq (-L)^v \circ \Sigma$. Then the same holds for their contact composition in the contactizations (cf.~\cite[Lemma 9.1.8]{LiNadlerShende}). Hence $(-L_{\cup,\epsilon}^-)^v \circ \Sigma_{\cup,\epsilon}^- \subseteq \Gamma_{\cup,\epsilon}^+ \cup (\Gamma_{\cup,\epsilon}^+)_{2\epsilon}$.\footnote{The second piece with $2\epsilon$ comes from the fact that composing two doublings results in four copies: if $(-L)^v \circ \Sigma = \Gamma$, then the contact composition satisfies $(-L_{\epsilon})^v \circ \Sigma_{-\epsilon'} = \Gamma_{\epsilon + \epsilon'}$.} 
We also know that, by a small negative push-off in the contactization, we may always replace $\Gamma = (-L)^v \circ \Sigma$ 
by its negative Reeb push-off and assume that $L \cap \mathfrak{c}_{M^- \times N} = \varnothing$ (cf.~\cite[Lemma 9.1.10]{LiNadlerShende}). Since $w_L^-(\cK)$ is also constructible with perfect stalks, $\Hom^\circ(w_{L}^-(\cK) \circ w_{\Lambda}^-(-))$ preserves colimits by Theorem \ref{thm:colim-preserve-hom}.

Lastly, we check that squeezability on Weinstein manifolds can be extended naturally along various embeddings, so that the main result in Section \ref{ssec: squeezability} applies:
Let $\Psi$ be the contact isotopy on $\partial^\infty M$ that squeezes $\Gamma^\infty$ into $\Lambda^\infty$. Let $M$ be a compact Weinstein domain whose completion gives the manifold. We can extend the contact isotopy $\Psi$ on $\partial^\infty M$ to a Hamiltonian isotopy in the domain $M$, and hence to a contact isotopy $\Psi_{M\times \real}$ in $M \times \real$. Consider coordinates $(x, r) \in \partial^\infty M \times \real_+$, and let $\rho: \real_+ \to [0, 1]$ be a bump function that equals $0$ on $(0, 1]$. We can set
\[H_{\Psi_{M\times \real}}(x, r, z; t) = r\rho(tr)H_{\Psi}(x),\]
which remains a negative isotopy with bounded energy $C$. It is straightforward to compute that if $\Psi$ squeezes $\Gamma^\infty$ to $\Sigma^\infty$, then $\Psi_{M\times \real}$ squeezes $\Gamma \times 0$ to $\Sigma \times C$.\footnote{While $\msh_{M}$ is defined as the subsheaf of $\msh_{M \times \real}$ consisting of microsheaves supported on $M \times \{0\}$, the shift to $C$ is simply a feature of isotoping microsheaves, and one can canonically identify the image back with $\msh_M$.}
Then we can extend it to an ambient contact isotopy $\Psi_{S^*\real^N}$ on $S^*\real^N$ by pulling back along the normal bundle and cutting off by a bump function, which remains a negative isotopy with bounded energy. Furthermore, it is straightforward from the construction that $\Gamma_{\cup,\epsilon}^+ \cup (\Gamma_{\cup,\epsilon}^+)_{2\epsilon}$ is squeezed into $\Lambda_{\cup,\epsilon}^+ \cup (\Lambda_{\cup,\epsilon}^+)_{2\epsilon}$ by ${\Psi}_{S^*\real^N}$. 
Now consider the diagram that factorizes $m_\Lambda \circ \iota_{\Lambda_{\cup,\epsilon}^+}^!$ as the composition:
\[\Shv_{\Gamma_{\cup,\epsilon}^+ \cup (\Gamma_{\cup,\epsilon}^+)_{2\epsilon}}(\real^N) \rightarrow \Shv_{\Lambda_{\cup,\epsilon}^+ \cup (\Lambda_{\cup,\epsilon}^+)_{2\epsilon}}(\real^N) \to \Shv_{\Lambda_{\cup,\epsilon}^+}(\real^N) \to \msh_\Lambda(\Lambda).\]
The first functor commutes with colimits by Proposition \ref{prop: squeezable-iota-shriek}. Since the second functor can be realized by the nearby cycle that takes the additional microsupport $\Lambda_{4\epsilon}$ to $\Lambda_{2\epsilon}$ by negative push-off, it preserves the microlocalization $m_\Lambda$ along $\Lambda$. Hence, the second and third functors agrees with $m_\Lambda$, and thus it also commutes with colimits by Proposition \ref{prop:musheaf-lim-colim}. 
This completes the proof. 
\end{proof}

\section{Cofinal wrapping dominated by the target} 

In this section, we will provide a purely geometric argument on Theorem \ref{thm:wrap-one-factor} in a special case. Namely, for a constructible sheaf $\cK \in \Shv_{\mathbb{R}-c}^b(X \times Y)$ with perfect stalks, suppose
$${\ms}^\infty(\cK) \subseteq T^*X \times S^*Y \setminus \Sigma^\infty.$$
Then for any constructible sheaf with perfect stalks $\cF \in \Shv_{\mathbb{R}-c}^b(X)$, the integral transform can be computed by the formula
$$\wrap_{-\Lambda \times \Sigma}^+\cK \circ \wrap_\Lambda^+\cF = \wrap_\Sigma^+(\cK \circ \cF),$$
by geometrically realizing the cofinal wrapping.


\subsection{Construction of cofinal wrappings}


Let us first make some observations. Let $M, N$ be Liouville manifolds. Denote by $\mathfrak{c}_M, \mathfrak{c}_N$ their cores (the attracting locus of the Liouville flow) and $\partial_\infty M, \partial_\infty N$ their ideal contact boundaries.
The boundary at infinity for the product $M\times N$ is $$\del^\infty(M\times N)=\del^\infty M\times \mathfrak{c}_N\cup \del^\infty M \times \del^\infty N \times \real \cup \mathfrak{c}_M \times \del^\infty N $$ with the contact form $\lambda_M + \lambda_N|_{\del^\infty(M\times N)}$. 

The following lemma implies that the condition of Gao in Theorem \ref{conj1} is a special case of our condition in the main Theorem \ref{thm:main}.

\begin{lem}\label{pre-flowon1factor}     
    Let $L \subseteq M \times N$ be a proper subset that is cylindrical at infinity. Then $L \to N$ is proper %
    if and only if $L^\infty \cap \partial^\infty M \times \mathfrak{c}_N = \varnothing$.
\end{lem}
\begin{proof}
    Properness of the map $L \to N$ %
    requires the cylindrical end $L^\infty \cap \partial^\infty M \times \mathfrak{c}_N = \varnothing$. Otherwise, we can find a compact subset $\overline{N}$ in $N$ such that $L \cap \pi_2^{-1}(\overline{N}) = L \cap M \times \overline{N}$ 
    contains $(L^\infty \cap \partial^\infty M \times \mathfrak{c}_N) \times \real_+$ which cannot be compact.

    Conversely, the empty intersection $L^\infty \cap \partial^\infty M \times \mathfrak{c}_N = \varnothing$ implies properness of the map $L \to N$. Otherwise, we can find a compact subset $\overline{N}$ in $N$ such that $L \cap \pi_2^{-1}(\overline{N}) = L \cap M \times \overline{N}$ is noncompact. Consider an exhaustion of $M$ by compact domains $M_n$ for $n \in \mathbb{N}$. Since $L$ is properly embedded, $L \cap M \times \overline{N}$ is noncompact means there is a sequence $(x_n, y_n) \in L$ such that $x_n \in M \setminus M_n$ and $y_n \in \overline{N}$. Consider the negative Liouville flow in $M \times N$. When $L$ is cylindrical outside a compact set $\overline{M} \times \overline{N}$, there exists $r_n > 0$ such that $(\varphi_Z^{-r_n}(x_n), \varphi_Z^{-r_n}(y_n)) \in L \cap \overline{M} \times \overline{N}$, where $\varphi_Z^{-r_n}(y_n) \in \varphi_Z^{-r_n}(\overline{N})$. Since $x_n \in M \setminus M_n$. The requirement $\varphi_Z^{-r_n}(x_n) \subseteq \overline{M}$ forces $r_n \to \infty$. Hence a subsequence $\varphi_Z^{-r_n}(y_n)$ converges to $\bigcap_{n=1}^\infty \varphi_Z^{-r_n}(\overline{N})$ which is equal to $\mathfrak{c}_N$. Since $L$ is proper, the limit is contained in $L^\infty \cap \partial^\infty M \times \mathfrak{c}_N$, which is a contradiction.
\end{proof}



\begin{construction}[The cofinal positive Hamiltonian flow]\label{flowon1factor}
Consider the stop 
$\mathfrak{c}_M \times \Sigma^\infty \cup \Lambda^\infty \times \Sigma^\infty \times \real \cup \Lambda^\infty \times \mathfrak{c}_N$
in the contact manifold $\partial^\infty(M \times N)$. We construct a cofinal positive Hamiltonian flow, 
which will be named $\tilde{H_2}$, away from the stop for any subset 
$$L^\infty \subseteq M \times \partial^\infty N \setminus \Sigma^\infty.$$
First, consider a cofinal Hamiltonian flow supported in $\partial^\infty N \setminus \Sigma^\infty$ and name it $H_2$:
\[H_2: \partial^\infty N \setminus \Sigma^\infty \to \real_{\ge 0}\] %
We use the coordinates $(x, r) \in \partial^\infty M \times \real_{>0}$ and $(y, s) \in \partial^\infty N \times \real_{>0}$, and set $r = 0$ on $\mathfrak{c}_M$ and $s = 0$ on $\mathfrak{c}_N$. Now we extend this flow to a Hamiltonian flow on $\partial^\infty(M \times N)$ satisfying the condition (which is well defined on $\mathfrak{c}_M \times \partial^\infty N$ where $r = 0$) 
\[\tilde{H}_2(x, r, y, s)=sH_2(y), \quad \text{when }\; (x,r, y, s) \in M \times (\partial^\infty N \times \real_{\geq \epsilon}).\] 
Let $U(\partial^\infty M \times \mathfrak{c}_N)$ be a neighborhood of $\partial^\infty M \times \mathfrak{c}_N \subseteq \partial^\infty(M \times N)$ such that
$$L^\infty \cap U(\partial^\infty M \times \mathfrak{c}_N) = \varnothing.$$
Let $\rho$ be a bump function on $\partial^\infty(M \times N)$ supported in $U(\partial^\infty M \times \mathfrak{c}_N)$ of a single variable $s/r \in [0, +\infty]$. We finally define
\[\tilde{H}_{2}(x, r, y, s) = (1 - \rho(s/r))H_2(y)\]
This is a homogeneous positive Hamiltonian on $\partial^\infty(M \times N) \times \real_{>0}$ and hence defines a contact Hamiltonian on $\partial^\infty(M \times N) \setminus M \times \Sigma^\infty$.
\end{construction}
\begin{rmk} 
Alternatively, one can also consider One can consider a bump function $\sigma: \real_{\geq 0} \to [0, 1]$ such that $\sigma(0) = 0$ and $\sigma(s) = s$ for $s \in \real_{\geq\epsilon}$, and then define the following Hamiltonian %
\[\tilde{H_2}(x, r, y, s)=(\sigma(s)^2{H}^2_{2}(y) + \sigma(r)^2\rho(s/r))^{1/2}\]
and verify that it satisfies all the properties we require. However, verifying that it is homogeneous outside a compact subset will be more delicate.
\end{rmk}

The second technical condition on the support of the bump function will make sure the image of $L^\infty$ is contained in the region where $\tilde{H}_2$ is exactly determined by $H_2$ and never enters the region where the bump function is non-zero:
\begin{lem}\label{steadyflow} Let ${U} = {U}(\partial^\infty M \times \mathfrak{c}_N) \subseteq \partial^\infty(M \times N)$ and let $L^\infty \subseteq \partial^\infty(M \times N)\setminus \overline{U}$ be any closed subset. Then for any $t \in \real$, we have $\varphi_{\tilde{H_2}}^t(L^\infty) \subseteq \partial^\infty(M \times N)\setminus \overline{U}$. 
\end{lem} \begin{proof}
We argue as follows. 
Let $z(t)$ be the Hamiltonian trajectory starting from any point $z_0 \in L^\infty$. Let $t^* = \inf \{ t \ge 0 \mid z(t) \in \overline{U} \} \in [0, +\infty]$ be the time for $z(t)$ to enter the region $\overline{U}$. 
For the entire interval $t \in [0, t^*)$, the trajectory $z(t)$ remains outside $\overline{U}$ where the bump function $\rho(s/r) = 0$, so the dynamics are governed by the Hamiltonian $\tilde{H}_2(x, r, y, s) = sH_2(y)$. We can rescale the contact form on $\partial^\infty N \setminus \Sigma^\infty$ such that we may assume after rescaling the Hamiltonian flow is the Reeb flow $H_2 = 1$. The equation of $z(t) = (x(t), r(t), y(t), s(t))$ satisfies the equations $r'(t) = 0$ and $s'(t) = 0$, which implies that $r(t) = r_0$ and $s(t) = s_0$ is constant for all $t \le t^*$. Consequently, we know that $s(t^*)/r(t^*) = s_0/r_0$. This shows $s(t^*)/r(t^*) \geq s_0/r_0$. Hence $z(t^*) \notin \overline{U}$ if $t^* < \infty$. This implies that $t^*=\infty$ and hence the trajectory is disjoint from $\overline{U}$.
\end{proof}

We then obtain the desired globally defined Hamiltonian flow $\tilde{H_2}$ over the product space $\partial^\infty(M \times N) \times \real_{>0}$, which outside of $U$ recovers the wrapping of $H_2$ on the target component $\partial^\infty N \times \real_{>0}$ and remains identity wrapping on the domain component $M$. %
Now we show the cofinality property for the flow we constructed. 
\begin{thm}\label{cofinal} For any closed subset $L^\infty \subseteq M \times \partial^\infty N \setminus \Sigma^\infty$, %
the positive Hamiltonian $\tilde{H_2}$ in Construction \ref{flowon1factor} generates a cofinal wrapping of $L^\infty \subseteq \partial^\infty(M \times N)$ in the wrapping category Definition \ref{def:wrapping-cat} with respect to the stop $\mathfrak{c}_M \times \Sigma^\infty \cup \Lambda^\infty \times \Sigma^\infty \times \real \cup \Lambda^\infty \times \mathfrak{c}_N$.
\end{thm} 
\begin{proof}
    It suffices to check the condition from \cite[Lemma 3.29]{GPS1} to show the cofinality. 
    Following \cite[Lemma 2.21]{GPS1}, pick a contact form on $\partial^\infty(M \times N) \setminus (\mathfrak{c}_M \times \Sigma^\infty \cup \Lambda^\infty \times \Sigma^\infty \times \real \cup \Lambda^\infty \times \mathfrak{c}_N)$ of the form %
    \begin{equation*}
    \alpha_{\partial^\infty(M \times N) } = \lambda_{\partial^\infty (M\times N)}+\lambda_{N}|_{\partial^\infty(M \times N) }.
    \end{equation*}
    Then, on the subset ${M} \times \partial^\infty N \setminus \Sigma^\infty$, we can write
    \begin{equation*}
    \alpha_{M \times \partial^\infty N \setminus \Sigma^\infty} = (r^2+s^2)^{-1/2}(\lambda_{M}+ s\alpha_{\partial^\infty N \setminus \Sigma^\infty}).
    \end{equation*}
    In Construction \ref{flowon1factor}, we can assume that $r/s \leq \epsilon$ outside of the support of the bump function $\rho$. By Lemma \ref{steadyflow}, since $L^\infty$ is outside of the support of $\rho$, we know $\varphi_{\tilde{H}_2}^t(L^\infty)$ is also outside of the support of $\rho$ for any $t \geq 0$. Without loss of generality, we may rescale the contact form $\alpha_{\partial^\infty N \setminus \Sigma^\infty}$ such that ${H}_2$ defines the Reeb flow on $\partial^\infty N \setminus \Sigma^\infty$. Then
    \begin{align*}
    \sup\alpha(X_{\tilde{H}_2})|_{\varphi_{\tilde{H}_2}^t(L^\infty)}
    &=\sup (r^2+s^2)^{-1/2}(\lambda_{M}(X_{\tilde{H}_2})+ s\alpha_{\partial^\infty N \setminus \Sigma^\infty}(X_{\tilde{H}_2}))|_{\varphi_{\tilde{H}_2}^t(L^\infty)}\\
    &=\sup (r^2+s^2)^{-1/2}s\alpha_{\partial^\infty N \setminus \Sigma^\infty}(X_{{H}_2})|_{\varphi_{\tilde{H}_2}^t(L^\infty)} \\
    &= (1 + (r/s)^2)^{-1/2} \geq (1 + 1/\epsilon^2)^{-1/2}.
    \end{align*} 
    Hence we have finished checking the condition in \cite[Lemma 3.29]{GPS1}.
\end{proof}

\subsection{Geometric proof of the theorem in a special case}\label{sec:special-case}
Using the geometric Construction \ref{flowon1factor} in the previous subsection, we give a new proof of Theorem \ref{thm:wrap-one-factor} in the special case when
$$\ms^\infty(\cK) \subseteq T^*X \times S^*Y \setminus \Sigma^\infty.$$

\begin{lem}\label{0preserveflow} 
Let $\Lambda \subseteq T^*X, \Sigma \subseteq T^*Y$ be closed subsets and $\cK \in \Shv(X \times Y)$ such that
\[\ms^\infty(\cK) \subseteq T^*X \times S^*Y \setminus \Sigma^\infty.\]
Then for any $\cF \in \Shv(X)$ and $\cG\in \Shv_\Sigma(Y)$, we have %
    \begin{align*}
        \hom(\wrap_{-\Lambda \times \Sigma}^+\cK,\Hom^\boxtimes(\cF, \cG)) \cong \hom( \cK,\Hom^\boxtimes(\cF, \cG)).  
    \end{align*} 
\end{lem}
    \begin{proof}
    Observe that the Hamiltonian flow on the product is eventually constant on the internal sheaf hom, in the sense that $$\cK(\phi_{\tilde{H}_{2}})|_t \circ \Hom^\boxtimes(\cF, \cG) = \Hom^\boxtimes(\cF, \cG).$$
    This is because we can estimate the microsupport by \cite[Proposition 5.4.2]{KS13}
    $$\ms(\Hom^\boxtimes(\cF, \cG))\subseteq -\ms(\cF) \times \Sigma.$$ %
    The cofinal Hamiltonian flow in Construction \ref{flowon1factor} acts as the identity map on $-\ms(\cF) \times \Sigma^\infty$. Consequently, we can invoke the constancy lemma \cite[Corollary 1.6]{GKS} or the uniqueness statement in Theorem \ref{thm:gks}, %
    which tells us the isomorphism. %
    Therefore, 
    \begin{align*}
        \hom(\cK,\Hom^\boxtimes(\cF, \cG))&\cong \operatornamewithlimits{colim}_{t \to \infty}  \hom(\cK,\Hom^\boxtimes(\cF, \cG))\\
        & \cong\operatornamewithlimits{colim}_{t \to \infty}\hom(\cK(\phi_{\tilde{H}_{2}})|_t\circ \cK,\Hom^\boxtimes(\cF, \cG)), 
        \end{align*}  
        where the second isomorphism holds due to the observation at the beginning of the proof.
    \end{proof}

\begin{proof}[Proof of Theorem \ref{thm:wrap-one-factor} in the special case]
With Lemma \ref{0preserveflow} now it follows that  %
\begin{align*}
\hom(\wrap_{-\Lambda \times \Sigma}^+\cK \circ \wrap_\Lambda^+\cF, \cG) &\cong \hom(\wrap_\Sigma^+(\wrap_{-\Lambda \times \Sigma}^+\cK \circ \cF), \cG) \\
&\cong \hom(\wrap_{-\Lambda \times \Sigma}^+\cK \circ \cF, \cG) \cong \hom(\wrap_{-\Lambda \times \Sigma}^+\cK, \Hom^\boxtimes(\cF, \cG)) \\
&\cong \hom(\cK, \Hom^\boxtimes(\cF, \cG)) \cong \hom(\cK \circ \cF, \cG) \\
&\cong \hom(\wrap_\Sigma^+(\cK \circ \cF), \cG),
\end{align*}
where the first isomorphism follows from Lemma \ref{lem: pushtoLambda} , the second and sixth isomorphisms follow from Theorem \ref{thm:wrap-adjoint}, and third and fifth isomorphisms follow from adjunctions of six functors, and the fourth isomorphism follows from Lemma \ref{0preserveflow}.
\end{proof}

Finally, we also remark that we can deduce Theorem \ref{thm:main} in the special case from Construction \ref{flowon1factor} plus Theorem \ref{thm:Hu} of Y.~Hu \cite{Hu25}.

\begin{proof}[Proof of Theorem \ref{thm:main} in the special case]
Let $\Lambda \subseteq T^*X, \Sigma \subseteq T^*Y$ be closed subsets and $\cK \in \Shv(X \times Y)$ such that
\[\ms^\infty(\cK) \subseteq T^*X \times S^*Y \setminus \Sigma^\infty.\]
First, we claim that for any $x \in X$, $\wrap_{-\Lambda \times \Sigma}^+\cK|_{\{x\} \times Y} = \wrap_\Sigma^+(\cK|_{\{x\} \times Y})$.
Since the cofinal Hamiltonian flow in Construction \ref{flowon1factor} acts as
\[p_{T^*Y}(\phi_{\tilde{H}_2}^t(\ms(\cK)) \cap T_x^*X \times T^*Y) = \phi_{H_2}(p_{T^*Y}(\ms(\cK) \cap T_x^*X \times T^*Y)).\]
Then, by the effect of the GKS kernel Proposition \ref{prop:gks-effect}, we have
\[(\cK(\phi_{\tilde{H}_2}) \circ \cK)|_{\{x\} \times Y \times I} = \cK(\phi_{H_2}) \circ (\cK|_{\{x\} \times Y}).\]
As we take the colimit for $t \to \infty$, this implies that
$$\wrap_{-\Lambda \times \Sigma}^+\cK|_{\{x\} \times Y} = \wrap_\Sigma^+(\cK|_{\{x\} \times Y}).$$
Since $\cK \in \Shv_{\real-c}(X \times Y)$, we know that $\cK|_{\{x\} \times Y} \in \Shv_{\real-c}(Y)$. Then by
Theorem \ref{wrap_compare}, we have
$$\wrap_{-\Lambda \times \Sigma}^+\cK|_{\{x\} \times Y} \cong \wrap_{\Sigma}^+(\cK|_{\{x\} \times Y})  \in \Shv^c_\Sigma(Y).$$
for any $x \in X$. Then the result follows from Theorem \ref{thm:Hu} of Y.~Hu \cite{Hu25}.
\end{proof}

\appendix

\section{Computational proof of Hu's result} \label{sec: Hu}

In this appendix, we give a computational proof of Hu's main result in \cite{Hu25}, stated in this paper as Theorem \ref{thm:Hu}. Hu's approach is an ingenious application of enriched category theory: instead of the usual approach of enhancing the coefficient from spectra to a more exotic coefficient $\cC$, he observes that one can view the target category $\Shv_{\Sigma}(Y)$ itself as a coefficient category $\cC'$, and then deduce the main theorem from general results.

The purpose of this appendix is to show that, for readers who are familiar with six-functor formalism, one can carry out the following computational proof. The authors do not assert that this is an alternative proof, since doing so would require tracing through the six-functor realization of the abstract argument; rather, we only try to present something more aligned to the spirit of \cite{KS13}. We will need a lemma:

\begin{lem} \label{lem: shriek-to-a-point}
Let $X$ be a real analytic manifold. For $x \in X$, the functor $i_x^!: \Shv(X) \rightarrow \Shv(\{x\}) = \cC$ restricts to a colimit-preserving functor on $\Shv_\Lambda(X)$ for any subanalytic $\Lambda \subseteq S^* X$ and is jointly conservative, i.e., for any $\cF \in \Shv_\Lambda(X)$, $i_x^! \cF = 0$ for all $x \in X$ if and only if $\cF = 0$.
\end{lem}

\begin{proof}
Let $x \in X$. Note that $\ms(1_x) = T^*_x X$ is a subanalytic Lagrangian; by \cite{Czapla}, one can always find a $C^1$ subanalytic Whitney triangulation $\cS$ containing both $\Lambda$ and $T^*_x X$. By Theorem \ref{thm: ss-constructible} and \cite[Lemma 4.5]{GPS3}, $1_x \in \Shv_{N^*(\cS)}(X)$ is compact\footnote{Instead of arguing using existence of nice triangulation, one could alternatively deduce it from the theory of exit-path spaces as in  \cite{Hu25}.} so $\iota_\Lambda^*(1_x) \in \Shv_\Lambda(X)$ is compact by Theorem \ref{thm: adjoint}. Since \[i_x^!(\cF) = \hom(1_{x}, \cF) = \hom(\iota_\Lambda^*(1_{x}), \cF)\] for $\cF \in \Shv_\Lambda(X)$, we see that $i_x^!$ is colimit-preserving. 

Now assume $\cF \in \Shv_\Lambda(X)$ satisfies $i_x^! \cF = 0$ for all $x \in X$, and we want to show that $\cF = 0$. By choosing a $C^1$ subanalytic Whitney triangulation $\cS$ such that $\Lambda \subseteq N^*(\cS)$, we can view $\cF$ as an object in $\Shv_{\cS}(X)$. Since we have already shown that $i_x^!$ is colimit-preserving and $\Shv_{\cS}(X)$ is compactly generated, it suffices to prove the case when $\cF$ is a compact object in $\Shv_{\cS}(X)$. Again by \cite[Lemma 4.5]{GPS3}, $\cF$ is in particular a constructible sheaf with perfect stalks (and compact support). 

The Verdier duality functor $D_X (-) \coloneqq \Hom(-,\omega_X)$ restricts to an equivalence on such objects. More precisely, if $\cF$ is a constructible sheaf with perfect stalks, then $D_X(\cF)$ is also constructible with perfect stalks. Furthermore, the double dual $D_X^2 \cF = \cF$ recovers the original sheaf. We then use \cite[Exercise VIII.3]{KS13} to compute that
\[i_x^! \cF = i_x^! D_X^2(\cF) = D_{\{x\}}( i_x^* D_X(\cF)).\]
The vanishing assumption then implies that $i_x^* D_X(\cF) = 0$ for all $x$ so $D_X(\cF) = 0$. Taking $D_X$ again, we see that $\cF = 0$.
\end{proof}

\begin{proof}[Computational proof of Theorem \ref{thm:Hu}]
We first show that if $\cK \in \Shv_{-\Lambda \times \Sigma}(X \times Y)$ induces a functor preserving compact objects, then $\cK|_{\{x\} \times Y}$ is compact in $\Shv_\Sigma^c(Y)$ for all $x \in X$. First, by our assumption, $\cK \circ \iota^*(1_x)$ is compact, but we can compute that
\[\cK \circ \iota^*(1_x) = \cK \circ 1_x = (\pi_Y)_! (\cK \otimes 1_{ \{x\} \times Y}) = \cK_{\{x\} \times Y}\]
as desired. Here, we use Lemma \ref{lem: niceformwrapping} for the first equality.
We now consider the other implication. As in the proof of Theorem \ref{thm:main}, we show that $\Hom^\circ(\cK,-)$ is colimit-preserving, i.e., for any family of sheaves $\cG_\alpha \in \Shv_{\Sigma}(Y)$, the canonical map
\begin{equation} \label{for: comparison-map-Hu}
\oplus_\alpha \Hom^\circ(\cK, \cG_\alpha) \rightarrow  \Hom^\circ(\cK, \oplus_\alpha \cG_\alpha)
\end{equation}
is an equivalence. We begin with the computation,
\[ i_x^! \Hom^\circ(\cK, \cG) = i_x^! {\pi_X}_* \Hom(\cK, \pi_Y^! \cG) = {a_X}_* i_x^{'!} \Hom(\cK, \pi_Y^! \cG) = \hom(\cK |_{\{x\} \times Y}, \cG).\]
Here, we use base change for the pullback diagram
\[
\begin{tikzcd}
	\{x\} \times Y && X \times Y  \\
	\{x\} && X
	\arrow["i_x' ", from=1-1, to=1-3]
	\arrow["i_x", from=2-1, to=2-3]
	\arrow["a_X ", from=1-1, to=2-1]
	\arrow["\pi_X", from=1-3, to=2-3]
\end{tikzcd}
\]
and the fact that $f^! \Hom(-,-) = \Hom( f^*(-), f^!(-))$.  By the previous Lemma \ref{lem: shriek-to-a-point}, we can check the equivalence \eqref{for: comparison-map-Hu} after applying $i_x^!$, since a map is an equivalence if and only if its cofiber vanishes. Since $i_x^!$ commutes with direct sums, this computation thus implies that it suffices to check that the map
\[\oplus_\alpha \hom^\circ(\cK |_{\{x\} \times Y} , \cG_\alpha) \rightarrow  \Hom^\circ(\cK |_{\{x\} \times Y}, \oplus_\alpha \cG_\alpha)\]
is an equivalence, which is exactly the assumption that $\cK |_{\{x\} \times Y}$ is compact.
\end{proof}


\bibliographystyle{amsalpha}
\bibliography{ref}

\end{document}